\documentclass[12pt]{article}
\usepackage[a4paper,dvips,twoside]{geometry}
\usepackage[latin1]{inputenc}
\usepackage[T1]{fontenc}
\usepackage[all]{xy}
\usepackage{theorem}
\usepackage{accents}
\usepackage{amssymb}

\usepackage{fancyhdr}
\newcommand{\BB}[1]{\mathbb{#1}}
\newcommand{\FF}{\BB{F}}
\newcommand{\RR}{\BB{R}}
\newcommand{\NN}{\BB{N}}
\newcommand{\Cal}{\mathcal}
\newcommand\he { {\stackrel{\mathrm{h.e}}{\simeq} } }
\def\ibar {\underline{i}}
\long\def\pull#1 {\left | \vphantom { ^\bullet {#1} } \right .\mskip -8 mu {   \overline { ^\bullet #1 } } }
\newcommand\In[1]{\accentset{\circ}{#1}}

\theoremstyle{break}
\newtheorem{bigdef}{Définition}
\newtheorem{bigth}{Théorème}

\newtheorem{PROP}{Proposition}[section]
\newtheorem{THEO}[PROP]{Théorème}
\newtheorem{LEMM}[PROP]{Lemme}
\newtheorem{CORO}[PROP]{Corollaire}
\newtheorem{CONJ}[PROP]{Conjecture}
\newtheorem{DEFI}[PROP]{Définition}
\newtheorem{EXEM}[PROP]{Exemple}
\newtheorem{REMA}[PROP]{Remarque}

\long\def\Prop #1\par{\begin{PROP} #1 \end{PROP}\par}
\long\def\Th #1\par{\begin{THEO} #1 \end{THEO}\par}
\long\def\Lemme #1\par{\begin{LEMM} #1 \end{LEMM}\par}
\long\def\Cor #1\par{\begin{CORO} #1 \end{CORO}\par}
\long\def\Conj #1\par{\begin{CONJ} #1 \end{CONJ}\par}
\long\def\Def #1\par{\begin{DEFI} #1 \end{DEFI}\par}
\long\def\Exp #1\par{\begin{EXEM} #1 \end{EXEM}\par}
\long\def\Rem #1\par {\begin{REMA} #1 \end{REMA}\par}

\def\Dem {\bigbreak \noindent \bf Démonstration.\enspace \rm}
\long\def\Demr #1: {\bigbreak \noindent \bf Démonstration (de #1).\enspace \rm}
\def\qed{\hfill \vbox { \hrule \hbox{\vrule height 8pt  \hskip 8pt \vrule height 8pt} \hrule}\medbreak}

\begin{document}
\title{
Invariance homotopique \\
de certains espaces de configurations
}
\author{Jean-Philippe~Jourdan\\
\footnotesize Département de Mathématiques, UMR 8524, Université de Lille 1\\
\footnotesize 59655 Villeneuve d'Ascq Cedex, France\\
\footnotesize \texttt{jourdan@math.univ-lille1.fr}
}
\date{}
\maketitle
\begin{abstract}
Pour une variété lisse $A$, on considère $\FF_k(A\times\RR)$ l'espace des confi\-gu\-rations ordonnées de $k$ particules distinctes dans $A\times\RR$. On effectue une construction explicite de l'espace de configurations $\FF_k(A\times\RR)$ et de la  suspension $(k-2)$-ième de $\FF_k(A)$. Puis l'on montre que, sous certaines conditions, le type d'homotopie de ces deux espaces ne dépend que de celui de $A$.
\bigbreak

For a smooth manifold $A$, we consider the ordered configuration space $\FF_k(A\times\RR)$ of  $k$ distinct points in
$A\times\RR$. We obtain an explicit homotopy construction of the configuration space $\FF_k(A\times\RR)$ and of the $(k-2)$-fold suspension of $\FF_k(A)$. Under certain conditions, we then show that the homotopy types of these two spaces depend only on the homotopy type of $A$.
\end{abstract}
\section{Introduction}
\label{sec:Intro}

Pour une variété $M$, l'espace de configurations de $k$ particules dans $M$ est défini par :
  $$ \FF_k (M) = \{(x_1,\dots,x_k)\in M^{\times k} \mid i\neq j \Rightarrow x_i\neq x_j \}.$$
Cet espace n'est pas un invariant du type d'homotopie de $M$, i.e. si les variétés $M$ et $N$ ont même type d'homotopie, les espaces $\FF_k (M)$ et $\FF_k (N)$ n'ont pas toujours le même type d'homotopie.
Par exemple, l'espace euclidien de dimension $p$, $\RR^p$, a  toujours le type d'homotopie d'un point quelque soit $p$. Or l'espace des configurations de $2$ points dans $\RR^p$ a le type d'homotopie d'une sphère de dimension $p-1$ et dépend donc de $p$.

On dispose cependant de résultats positifs dans certains cas particuliers. Dans \cite{Levitt1995}, N.~Levitt montre que l'espace des configurations de $2$ particules, $\FF_2(M)$, d'une variété compacte et $2$-connexe $M$, est un invariant du type d'homotopie de $M$. Dans \cite{Lambrechts2001}, P.~Lambrechts et D.~Stanley  ont obtenu un résultat similaire en homotopie rationnelle, i.e. le type d'homotopie rationnelle de $\FF_2(M)$ ne dépend que du type d'homotopie rationnelle de $M$ si $M$ est une variété compacte et $2$-connexe. 
Toujours en homotopie rationnelle mais pour un nombre arbitraire de parti\-cu\-les, B.~Totaro \cite{Totaro1996} et I.~K{\v{r}}{\'{\i}}{\v{z}} \cite{Kriz1994} montrent que le type d'homotopie rationnelle de $\FF_k (M)$ ne dépend que de l'anneau de cohomologie de $M$ si $M$ est une variété projective complexe.

 Néanmoins, il n'existe pas de résultat général : dans \cite{Longoni2004}, R.~Longoni et P.~Salvatore exhibent deux espaces lenticulaires, $L_{7,1}$ et $L_{7,2}$, ayant même type d'homotopie mais dont les espaces de configurations ont un type d'homotopie différent. Le problème d'invariance du type d'homotopie de l'espace $\FF_k(M)$ dans le cas où $M$ est une variété compacte et simplement connexe (ou même $p$-connexe) reste cependant toujours ouvert.

Si l'on s'intéresse non pas à l'espace de configurations $\FF_k(M)$ lui-même, mais à certaines constructions effectuées à partir de $\FF_k(M)$, plusieurs résultats existent déjà dans la littérature :
\begin{itemize}
	\item $\Omega\FF_k (M)=\mathrm{Maps}_*(S^1,\FF_k (M))$ : \emph{l'espace des lacets pointés} de $\FF_k (M)$. Pour un espace pointé $X$, $\Omega X$ est défini comme l'espace des applications allant du cercle (pointé) $S^1$ dans $X$ et qui conservent le point de base .
 N.~Levitt \cite{Levitt1995} montre que le type d'homotopie de $\Omega\FF_k (M)$ ne dépend que de celui de $M$ si $M$ est une variété compacte. La structure de H-espace de $\Omega\FF_k (M)$ est détaillée par F.~R.~Cohen et S.~Gitler \cite{Cohen2002} pour des variétés particulières, notamment pour des variétés tressables («braidable manifolds»).
  	\item $\sum^{\alpha} \FF_k (M)$ : \emph{la suspension itérée} de $\FF_k (M)$. Pour un espace non pointé $X$, la suspension itérée $\sum^{\alpha} X$ est définie comme la somme amalgamée de $S^{\alpha-1}$ et $X \times D^\alpha$ le long de $X \times S^{\alpha-1}$. Dans le cas où $M$ est une variété compacte, J.~R.~Klein et M.~Aouina \cite{Aouina2003} montrent que $\sum^{\alpha} \FF_k (M)$ est un invariant du type d'homotopie de $M$ pour un $\alpha$ bien choisi. Si l'on suppose $M$ $2$-connexe et de dimension $d$, alors le  nombre de suspensions nécessaires est $\alpha=(k-2)d$.
\end{itemize}

 Dans la suite, on va s'intéresser à ces deux aspects.
Tout d'abord, on effectue une construction homotopique explicite des espaces de configurations pour une variété produit particulière $M=A\times\RR$.
Ce cas est l'exemple fondamental de variété tressable (voir \cite{Cohen2002}).
La construction de l'espace de configurations à $k$ particules $\FF_k(A\times\RR)$ a la particularité de s'effectuer à partir de $A$ et $\FF_2(A)$ pour tout $k$.
Ceci nous permet d'établir l'invariance homotopique de $\FF_k(A\times\RR)$ dans un cas particulier :
\begin{bigdef}
Considérons deux variétés $A$ et $B$ et une équivalence d'homotopie $f:A\to B$. L'application $f$ est appelée \emph{$\Cal F_2$-équivalence d'homotopie} s'il existe une équivalence d'homotopie $\varphi_f : \FF_2(A)\to \FF_2(B)$ faisant commuter à homotopie près le diagramme suivant :
 $$
 \xymatrix{
 {\FF_2 (A)}\ar[d]_{\varphi_f}\ar@{^(->}[r]	&	A\times A\ar[d]^{f\times f}	\\
 {\FF_2 (B)}\ar@{^(->}[r]		&	B\times B	\\
 }
 $$
\end{bigdef}
Remarquons d'après \cite{Aouina2003}, que toute équivalence d'homotopie entre deux variétés com\-pactes, $2$-connexes est une $\Cal F_2$-équivalence d'homotopie, cf. exemple \ref{F2equiv2}. À l'aide de notre construction de $\FF_k(A\times\RR)$, nous établissons alors :
\begin{bigth}
Soit $f:A\to B$ une $\Cal F_2$-équivalence d'homotopie entre deux variétés $A$ et $B$, alors pour tout $k\ge 1$, on a une équivalence d'homotopie fibrée :
$$\xymatrix{
{\FF_{k+1} }(A\times \RR)\ar[rr]^\he\ar@{->>}[d]	&& {\FF_{k+1} }(B\times \RR)\ar@{->>}[d] \\
{\FF_k }(A\times \RR)\ar[rr]_\he	&& {\FF_k }(B\times \RR) \\}
$$
\end{bigth}

Nous utilisons ensuite ce résultat ainsi que notre construction de $\FF_k(A\times\RR)$ pour étudier une suspension itérée des espaces de configurations de $A$~: $\sum^\alpha\FF_k(A)$.
Nous en  déduisons alors le théorème d'invariance homotopique suivant :
\begin{bigth}
Soit $f:A\to B$ une $\Cal F_2$-équivalence d'homotopie entre deux variétés connexes $A$ et $B$, alors pour tout $k\ge 1$, les espaces $\Sigma^{k-2} \FF_k (A)$ et $\Sigma^{k-2}\FF_k (B)$ ont même type d'homotopie.
\end{bigth}
Cet énoncé améliore le résultat de M.~Aouina et J.~R.~Klein \cite{Aouina2003}, l'ordre de la suspension itérée requise étant ici indépendant de la dimension des variétés intervenant.

Ce travail se déroule de la manière suivante : dans la section \ref{sec:Prelim}, nous donnons les définitions de base et présentons l'espace des configurations de $k$ points sur la droite réelle. Dans la section \ref{sec:Construct}, pour toute variété $A$, nous construisons un espace $E^k(A)$ de même type d'homotopie que l'espace des configurations de $k$ points dans $A\times \RR$. Les sections \ref{sec:Invariant} et \ref{sec:Suspension} sont consacrées respectivement aux preuves des théorèmes A et B. Finalement, un appendice regroupe les définitions et propriétés de limites homotopiques utilisées dans cette rédaction.

\section{Préliminaires}
\label{sec:Prelim}

 Afin de simplifier les démonstrations, nous allons introduire une variante des espaces de configurations. Cette définition se montrera bien adaptée aux démonstrations par récurrence que nous utiliserons dans cette section et les suivantes.

\Def Pour tout espace $A$ muni d'une filtration croissante $A_\bullet=(A_i)_{i\in\NN^*\cup\{\infty\}}$ avec $A_\infty=A$, nous définissons les espaces de configurations suivants :
$$\FF_k (A_\bullet)=\{ (a_1,\dots,a_k)\in A_1\times\dots A_{k-1}\times A_k \mid i\ne j \Rightarrow a_i\ne a_j \}$$
$$\overline{\FF}_k (A_\bullet)=\{ (a_1,\dots,a_k)\in A_1\times\dots\times A_{k-1}\times \overline{A_k} \mid i\ne j \Rightarrow a_i\ne a_j \}$$ où $\overline{A_k}$ est l'adhérence de $A_k$ dans $A$. Si de plus la filtration $A_\bullet$ vérifie $A_{*-1}\subseteq \In A_*$, on pose :
$$\In\FF_k (A_\bullet)=\{ (a_1,\dots,a_k)\in A_1\times\dots\times A_{k-1}\times \In A_k \mid i\ne j \Rightarrow a_i\ne a_j \}$$
et
$$\partial\FF_k (A_\bullet)=\overline{\FF}_k (A_\bullet)\setminus\In\FF_k (A_\bullet)=\FF_{k-1} (A_\bullet)\times\partial A_k .$$

Pour la filtration constante $A_\bullet$ définie par $A_i=A$, les trois espaces $\FF_k (A_\bullet)$, $\overline{\FF}_k(A_\bullet)$ et $\In\FF_k (A_\bullet)$ coïncident et nous retrouvons la définition usuelle des espaces de configurations de $A$. Nous noterons simplement $A$ cette filtration triviale.

\Def Si $A$ est une variété, on note $A\times J$ la filtration de $A\times \RR$ donnée par 
$$A\times\{0\}\subseteq A\times[-1,1]\subseteq\cdots\subseteq A\times[-k,k]\subseteq\cdots \textrm{ et } (A\times J)_\infty=A\times \RR .$$
Les espaces de configurations associés à cette filtration sont notés $\FF_k(A\times J)=\overline{\FF}_k(A\times J)$ et $\In\FF_k(A\times J)$.
Dans le cas où $A$ est réduite à un point $*$, nous noterons simplement $J$ la filtration de l'espace des nombres réels donnée par $*\times J$. 

Soit $A_\bullet$ une filtration associée à une variété $A$. Si $k$ et $r$ sont deux entiers tels que $r<k$, on note $\pi_{k,r}:\FF_k(A_\bullet)\to\FF_r(A_\bullet)$ la projection sur les $r$ premières composantes :
$$\pi_{k,r} (x_1,\dots,x_k)=(x_1,\dots,x_r).$$
On note de la même manière la projection $\pi_{k,r}:\overline{\FF}_k(A_\bullet)\to\FF_r(A_\bullet)$.

Dans  \cite{Fadell1962}, E.~Fadell et L.~Neuwirth montrent que si $A$ est une variété, alors $\pi_{k,r} : \FF_k(A)\to\FF_r(A)$ est un fibré localement trivial dont la fibre au dessus d'un élément $(q_1,\dots,q_r)$ de $\FF_r(A)$ est~:
 $$\FF_{k-r}(A\setminus\{q_1,\dots,q_r\}).$$
On a des résultats similaires pour les espaces de configurations associés à certaines filtrations, en particulier nous utiliserons les propriétés suivantes :
\Prop Soit $A$ une variété lisse.
\begin{enumerate}
	\item  Pour tout entier $k$, $k\ge 1$, la projection $\pi_{k+1,k} : \FF_{k+1}(A\times J)\to \FF_k(A\times J)$ est 
	un fibré localement trivial dont la fibre au dessus d'un élément $(q_1,\dots,q_k)$ de $\FF_k(A\times J)$ est :
 $$(A\times [-k,k])\setminus\{q_1,\dots,q_k\}.$$
	\item  Pour tout couple d'entiers $(k,r)$, $k>r\ge 0$, les inclusions $\FF_* (A\times J)\subset\FF_* (A\times\RR)$ induisent des équivalences d'homotopie fibrée :
	$$\xymatrix{
	{\FF_k (A\times J)}\ar@{^(->}[r]^\he\ar[d]		&{\FF_k (A\times \RR)}\ar[d]	\\
	{\FF_r (A\times J)}\ar@{^(->}[r]^\he		&{\FF_r (A\times \RR)}	\\
	}$$
\end{enumerate}

\Dem La première assertion est une adaptation immédiate de la dé\-monstra\-tion de E.~Fadell et L.~Neuwirth donnée dans \cite{Fadell1962}. Celle-ci permet de prouver  par récurrence que les inclusions $\FF_* (A\times J)\subset\FF_* (A\times\RR)$ sont des équivalences d'homotopie. La commutativité du carré de la proposition se vérifie aisément à partir de la définition des espaces $\FF_* (A\times J)$.\qed

En appliquant ce résultat à la filtration $J$ de l'espace des nombres réels $\RR$, on obtient une équivalence d'homotopie entre $\FF_k(J)$ et le produit $\FF_{k-1}(J)\times([-k+1,k-1]\setminus\{q_1,\dots,q_{k-1}\})$. Par récurrence, on constate que $\FF_k(J)$ a le type d'homotopie d'un espace discret à $k!$ points. Nous allons l'exprimer sous une forme mieux adaptée aux constructions à venir.
Pour $k\ge 2$, on pose :
$$\Cal A _k =\{1,3\}\times\{1,3,5\}\times\{1,3,5,7\}\times\cdots\{1,3,\dots,2k-3\}\times\{1,3,\dots,2k-1\}.$$
On pose aussi $\Cal A_1 = \emptyset$. On construit une équivalence d'homotopie entre $\FF_k (J)$ et $\Cal A_k$ par récurrence sur $k$.
Pour $k=2$, on définit $\theta_2 :\Cal A_2 \to \FF_2 (J)$ par $\theta_2(1)=(0,-1)$ et $\theta_2(3)=(0,+1)$.

Supposons avoir construit $\theta_k:\Cal A_k \to \FF_k (J)$, et choisissons $(i_2,\dots,i_k,2\alpha+1)$ un élément de $\Cal A_{k+1}$. Posons $\theta_k(i_2,\dots,i_k)=(t_1,\dots,t_k)$, et notons $[\![1,k]\!]$ l'ensemble des entiers naturels compris entre $1$ et $k$. Nous commençons par réordonner les éléments $(t_1,\dots,t_k)$ en choisissant des éléments $j_1,\dots,j_k$ dans $[\![1,k]\!]$ tels que $t_{j_1}<\cdots<t_{j_k}$. Notons alors $\sigma_k(i_2,\dots,i_{k})=(j_1,\dots,j_k)$.

Nous allons compléter $(t_1,\dots,t_k)$ par un élément $t_{k+1}$ qui va être choisi de la façon suivante :
\begin{itemize}
	\item Si $\alpha=0$, on pose $t_{k+1}=-k$.
	\item Si $\alpha=k$, on pose $t_{k+1}=+k$.
	\item Si $1\le \alpha \le k-1$, on choisit $t_{k+1}$ tel que $t_{j_{\alpha}}<t_{k+1}<t_{j_{\alpha+1}}$.
\end{itemize}

On pose alors : 
$$\theta_{k+1}(i_2,\dots,i_k,2\alpha+1)=(t_1,\dots,t_{k+1}).$$

Les applications $\theta_k$ sont clairement des équivalences d'homotopie qui rendent commutatif le carré suivant :
$$\xymatrix{
{\Cal A_{k+1}}\ar[rr]^{\theta_{k+1}}\ar[d]	&&	{\FF_{k+1}(J)}\ar[d]\\
{\Cal A_{k}}\ar[rr]^{\theta_{k}}		&&	{\FF_{k}(J)}\\
}$$

Avec les notations précédentes, on a défini une bijection $\sigma_k$ entre ${\Cal A_{k}}$ et le groupe des permutations $\Sigma_k$ par $\sigma_k(i_2,\dots,i_{k})(p)=j_p$ pour tout $p\in
[\![1,k]\!]$.

\section{Construction de $\FF_{k} (A\times\RR)$}
\label{sec:Construct}

Dans ce paragraphe, nous fixons une variété lisse $A$. Nous allons montrer qu'une construction homotopique de l'espace de configurations d'un nombre quelconque de particules dans $A\times\RR$ est toujours possible à partir des espaces $A$ et $\FF_2(A)$. 

Soit $\ibar=(i_2,\dots,i_{k-1})$ un élément de $\Cal A_{k-1}$, $k\ge 2$, et soit $(\ibar,i_k)=(i_2,\dots,i_{k-1},i_k)$ un élément de $\Cal A_{k}$.
Nous posons :
$$E^k_{i_2,\dots,i_k}=E^k_{\ibar,i_k}=A^k$$
et définissons une application :
$$\phi^k_{i_2,\dots,i_k}=\phi^k_{\ibar,i_k}:E^k_{\ibar}\to\FF_k (A\times J)$$
par $\phi^k_{\ibar,i_k}(x_1,\dots,x_k)=T_k(x_1,\dots,x_k,\theta_k(\ibar,i_k))$ où $(x_1,\dots,x_k,\theta_k(\ibar,i_k))\in A^k\times\FF_k(J)$ et $T_k : A^k\times\FF_k(J)\hookrightarrow\FF_k (A\times J)$, $T_k((a_1,\dots,a_k),(t_1,\dots,t_k))=((a_1,t_1),\dots,(a_k,t_k))$, est l'injection canonique.

Nous posons aussi $E^1=A$ et définissons $\phi^1:E^1\to\FF_1 (A\times J)$ par $\phi^1(a)=(a,0)$.
Remarquons que $\phi^1$ est une équivalence d'homotopie.

Considérons maintenant un élément $(\ibar,2p)\in\Cal A_{k-1}\times\{2,4,\dots ,2k-2\}$ et posons
$$E^k_{i_2,\dots,i_{k-1},2p}=E^k_{\ibar,2p}=\{(a_1,\dots,a_k)\in A^k \mid a_k \ne a_{j_p} \}$$
où l'élément $j_p$ est défini par $\sigma_{k-1}(i_2,\dots,i_{k-1})=(j_1,\dots,j_{k-1})$.
Pour rendre plus explicites les deux indices pour lesquels les coordonnées sont distinctes, nous notons aussi $E^k_{\ibar,2p}=A^k_{k j_p}$. Cet espace est évidemment homéomorphe à $\FF_2(A)\times A^{k-2}$.

Par construction, nous avons des inclusions naturelles $E^k_{\ibar,2p}\to E^k_{\ibar,2p-1}$ et $E^k_{\ibar,2p}\to E^k_{\ibar,2p+1}$.  Nous définissons alors $E^k_{\ibar}$, $\ibar\in\Cal A_{k-1}$,  comme la colimite homotopique de :
 $$ E^k_{\ibar,1}\gets E^k_{\ibar,2}\to E^k_{\ibar,3}\gets E^k_{\ibar,4}
 \gets\cdots\to
 E^k_{\ibar,2k-3}\gets E^k_{\ibar,2k-2}\to E^k_{\ibar,2k-1}.$$

\Prop \label{pull1} Soit $\ibar=(i_2,\dots,i_{k-1}) \in \Cal A_{k-1}$.
Les applications $\phi^k_{\ibar,*}:E^k_{\ibar,*}\to \FF_k (A\times J)$, définies ci-dessus pour tout entier $*\in \{1,3,\dots,2k-1\}$, induisent une application 
$$\phi^k_{\ibar}:E^k_{\ibar}\to \FF_k (A\times J)$$
 telle que le carré suivant soit cartésien pour tout $k\ge3$ :
$$\xymatrix{
E^k_{\ibar}\ar[rr]^{\phi^k_{\ibar}}\ar[d]	&& {\FF_k (A\times J)}\ar[d] \\
E^{k-1}_{\ibar}\ar[rr]^{\phi^{k-1}_{\ibar}}	&& {\FF_{k-1} (A\times J)} \\
}
$$

\Exp \label{F2AxR} Si $k=2$ alors $E^2_1=A^2$, $E^2_2=A^2_{21}=\FF_2 (A)$ et $E^2_3=A^2$.
Or l'application $\phi^1:E^1\to\FF_1 (A\times J)$ est une équivalence d'homotopie. La proposition ci-dessus montre que la colimite homotopique de $A^2\gets \FF_2 (A)\to A^2$ a le même type d'homotopie que  $\FF_2 (A\times J)$.

\Demr \ref{pull1}: Considérons le diagramme suivant :
$$\xymatrix@R=2ex{
		&E^k_{\ibar 1}= A^k\ar@/^2pc/[rrrdddd]|{\phi^k_{\ibar 1}}	&\\
	E^k_{\ibar 2}=A^k_{kj_1}\ar[ur]\ar[dr]\\
			&E^k_{\ibar 3}=A^k\ar@/^1pc/[rrrdd]|{\phi^k_{\ibar 3}}	&\\
	E^k_{\ibar 4}=A^k_{kj_2}\ar[ur]\ar[dr]\\
\vdots		&E^k_{\ibar 5}=A^k\ar[rrr]|{\phi^k_{\ibar 5}}	&&&{\FF_k (A\times J)}\\
	E^k_{\ibar (2k-4)}=A^k_{kj_{k-2}}\ar[dr]			&{\vdots}&&& 	\\
			&E^k_{\ibar (2k-3)} = A^k\ar@/_1pc/[rrruu]|{\phi^k_{\ibar (2k-3)}}	&\\
	E^k_{\ibar (2k-2)}=A^k_{kj_{k-1}}\ar[ur]\ar[dr]\\
			& E^k_{\ibar (2k-1)} = A^k\ar@/_2pc/[rrruuuu]|{\phi^k_{\ibar (2k-1)}}	&\\
}$$

En notant $\theta_k(\ibar,2p-1)=(t_1,\dots,t_{k-1},t^{(p)}_k)$  et $\theta_k(\ibar,2p+1)=(t_1,\dots,t_{k-1},t^{(p+1)}_k)$, on a par définition de $\theta_k$ : $$t^{(1)}_k<t_{j_1}<\dots<t_{j_{p-1}}<t^{(p)}_k<t_{j_{p}}<t^{(p+1)}_k<t_{j_{p+1}}<\dots<t_{j_{k-1}}<t^{(k)}_k.$$
Ce diagramme commute à homotopie près. L'homotopie pour le $p$-ième carré, $1\le p \le k-1$, est donnée explicitement par :
$$\xymatrix@=1ex{
E^k_{\ibar,2p}\times[0,1]=A^k_{k j_p}\times[0,1]\ar[r]	&	{\FF_k (A\times J)}\\
(a_1,\dots,a_k,s)\ar@{|->}[r]	& ((a_1,t_1),\dots,(a_{k-1},t_{k-1}),(a_k,s t_k^{(p)}+(1-s)t_k^{(p+1)})\\
}$$
Cette application est bien définie car par définition de $A^k_{k j_p}$, le seul entier $\alpha\in[\![1,k-1]\!]$ tel que $t_k^{(p)}\le t_\alpha\le t_k^{(p+1)}$ est $\alpha=j_p$ et l'on a $a_k\ne a_{j_p}$ .
Par la propriété universelle d'une colimite homotopique, on obtient une application $\phi^k_{\ibar} : E^k_{\ibar} \to {\FF_k (A\times J)}$.

De plus, chacun des facteurs de cette colimite se projette sur $E^{k-1}_{\ibar}=A^{k-1}$ en gardant les $k-1$ premières coordonnées de manière à ce que le diagramme
$$\xymatrix{
E^k_{\ibar(2p - 1)}\ar[d] & E^k_{\ibar 2p}\ar[r]\ar[d]\ar[l]& E^k_{\ibar(2p + 1)}\ar[d] \\
E^{k-1}_{\ibar}          & E^{k-1}_{\ibar}\ar@{=}[r]\ar@{=}[l]     & E^{k-1}_{\ibar}
}$$
commute. Nous en déduisons la commutativité à homotopie près du carré de la proposition :
$$\xymatrix{
E^k_{\ibar}\ar[rr]^{\phi^k_{\ibar}}\ar[d]	&& {\FF_k (A\times J)}\ar[d] \\
E^{k-1}_{\ibar}\ar[rr]^{\phi^{k-1}_{\ibar}}	&& {\FF_{k-1} (A\times J)} \\
}
$$
Soit $(q_1,\dots,q_{k-1})\in E^{k-1}_{\ibar}=A^k$ d'image ${\phi^{k-1}_{\ibar}}(q_1,\dots,q_{k-1})=((q_1,t_1),\dots,(q_{k-1},t_{k-1}))$.
La fibre de ${\FF_k (A\times J)}\to{\FF_{k-1} (A\times J)}$ au dessus de ${\phi^{k-1}_{\ibar}}(q_1,\dots,q_{k-1})$ est
$$ A\times [-k+1,k-1]\setminus\{(q_1,t_1),\dots,(q_{k-1},t_{k-1})\}.$$
Rappelons que par définition de $\theta_k$, $\theta_k(\ibar,2p-1)=(t_1,\dots,t_{k-1},t^{(p)}_k)$.
La proposition \ref{Puppe2} nous montre que la fibre homotopique de $E^k_{\ibar} \to E^{k-1}_{\ibar}$ au dessus de $(q_1,\dots,q_{k-1})$ est la colimite homotopique du diagramme :
$$ A\gets A\setminus q_1\to A
 \gets\cdots\to
 A\gets A\setminus q_{k-1}\to A$$
c'est-à-dire l'espace 
$$A\times\{t^{(1)}_k,\dots,t^{(k)}_k\}\bigcup_{i=1\dots k-1}(A\setminus q_i)\times[t^{(p)}_k,t^{(p+1)}_k].$$
On vérifie facilement que la restriction de $\phi^k_{\ibar}$ aux fibres est une équivalence d'homotopie.
Les deux flèches verticales du carré de l'énoncé ayant même fibre homotopique, ce carré est cartésien (voir \cite[théorème 2.6]{May1975}).
\qed

\emph{Pour tout $k\ge 2$, tout $\ibar\in\Cal A_{k-1}\cup\Cal A_{k}$, nous disposons d'applications 
$$\phi^k_{\ibar} : E^k_{\ibar} \to \FF_k (A\times J)$$
Nous allons maintenant procéder à une récurrence pour définir 
$$\phi^k_{\ibar} : E^k_{\ibar} \to \FF_k (A\times J) \textrm{ pour } \ibar\in\Cal A_\alpha \textrm{ avec } 1\le \alpha \le k .$$}

Fixons $(i_2,\dots,i_\alpha)\in\Cal A_\alpha$. Pour tout $k$, $k\ge\alpha+1$, nous construisons, par récurrence sur le nombre $k-\alpha$, deux familles d'espaces :
$$(E^k_{i_2,\dots,i_\alpha,2p})_{p\in[\![1,\alpha]\!]} \textrm{ et } (E^k_{i_2,\dots,i_\alpha,2p-1})_{p\in[\![1,\alpha+1]\!]}.$$

Observons que l'espace $E^k_{i_2,\dots,i_{\alpha}}$ est de la forme $E^k_{i_2,\dots,i_{\alpha-1},2p-1}$ avec $(i_2,\dots,i_{\alpha-1})\in\Cal A_{\alpha-1}$ et $p\in[\![1,\alpha]\!]$. Nous définissons donc :
\begin{itemize}
	\item l'espace $E^k_{i_2,\dots,i_\alpha}$ comme la colimite homotopique de
$$\scriptstyle E^{k}_{i_2,\dots,i_\alpha,1}\gets E^{k}_{i_2,\dots,i_\alpha,2}\to E^{k}_{i_2,\dots,i_\alpha,3}
 \gets\cdots\to
 E^{k}_{i_2,\dots,i_\alpha,2\alpha-1}\gets E^{k}_{i_2,\dots,i_\alpha,2\alpha}\to E^{k}_{i_2,\dots,i_\alpha,2\alpha+1}$$	
	\item l'espace $E^{k+1}_{i_2,\dots,i_\alpha,2p}$ comme le produit fibré homotopique
$$\xymatrix{
E^{k+1}_{i_2,\dots, i_\alpha, 2p}\ar[rr]\ar[d]	&& E^{k+1}_{i_2,\dots, i_\alpha,2p-1}\ar[d] \\
E^{k}_{i_2,\dots ,i_\alpha, 2p}\ar[rr]		&& E^{k}_{i_2,\dots ,i_\alpha,2p-1} \\
}
$$
\end{itemize}
Dans la preuve ci-dessous, nous justifions la validité de ce procédé en exhibant les applications que nécessite la construction et démontrons :
 
\Prop \label{pull} Pour tout $(i_2,\dots,i_\alpha)\in\Cal A_\alpha$ et tout entier $k$, $k\ge\alpha+1$, il existe un carré cartésien
$$\xymatrix{
E^{k+1}_{i_2,\dots,i_\alpha}\ar[rr]^{\phi^{k+1}_{i_2,\dots,i_\alpha}}\ar[d]	&& {\FF_{k+1} (A\times J)}\ar[d] \\
E^{k}_{i_2,\dots,i_\alpha}\ar[rr]^{\phi^{k}_{i_2,\dots,i_\alpha}}	&& {\FF_k (A\times J)} \\
}
$$

\Dem Nous effectuons une démonstration par récurrence sur l'entier $ k - \alpha $.
Soit $s$ un entier positif. L'hypothèse de récurrence est la propriété $H_s$ suivante :
pour tout couple d'entiers $(k,\alpha)$ tel que $0\le k-\alpha\le s$, pour tout $(i_2,\dots,i_\alpha)\in\Cal A_\alpha$, et pour tout $p\in[\![1,\alpha]\!]$, on a construit les espaces  $E^k_{i_2,\dots,i_\alpha}$ et $E^{k+1}_{i_2,\dots,i_\alpha,2p}$ suivant le procédé décrit avant l'énoncé et des applications $\phi^{k+1}_{i_2,\dots,i_\alpha}$, $\phi^k_{i_2,\dots,i_\alpha}$ de sorte que le carré suivant soit cartésien :
$$\xymatrix{
	E^{k+1}_{i_2,\dots,i_\alpha}\ar[rr]^{\phi^{k+1}_{i_2,\dots,i_\alpha}}\ar[d]	&& {\FF_{k+1} (A\times J)}\ar[d] \\
	E^{k}_{i_2,\dots,i_\alpha}\ar[rr]^{\phi^{k}_{i_2,\dots,i_\alpha}}	&& {\FF_k (A\times J)} \\
	}
	$$

La proposition \ref{pull1} montre que la propriété $H_0$ est vraie.
Supposons maintenant vraie la propriété $H_{s-1}$ .

 Fixons donc un entier $\alpha$, $\alpha\ge 1$, et un élément $(i_2,\dots,i_\alpha)\in\Cal A_\alpha$. Posons $k=\alpha+s$. De la propriété $H_{s-1}$, on déduit que pour tout entier $p$, $1\le p\le \alpha$, le diagramme en traits pleins suivant est commutatif  :
$$\xymatrix{
E^{k+1}_{i_2,\dots ,i_\alpha, 2p}\ar[rr]\ar[dd]\ar@{..>}[rd]	&& E^{k+1}_{i_2,\dots ,i_\alpha,(2p-1)}\ar[rr]^{\phi^{k+1}_{i_2,\dots ,i_\alpha,(2p-1)}}\ar[dd] && {\FF_{k+1} (A\times J)}\ar[dd]\\
&E^{k+1}_{i_2,\dots ,i_\alpha,(2p+1)}\ar[urrr]_>(0.7){\phi^{k+1}_{i_2,\dots ,i_\alpha,(2p+1)}}\ar[dd]\\
E^{k}_{i_2,\dots ,i_\alpha ,2p}\ar[rr]\ar[dr]		&& E^{k}_{i_2,\dots ,i_\alpha,(2p-1)}\ar[rr]^{\phi^{k}_{i_2,\dots ,i_\alpha,(2p-1)}} && {\FF_k (A\times J)}\\
&E^{k}_{i_2,\dots ,i_\alpha,(2p+1)}\ar[urrr]_{\phi^{k}_{i_2,\dots ,i_\alpha,(2p+1)}}\\
}
$$
et que chaque carré est cartésien. D'après le lemme du prisme (\ref{lemprism}), il existe une application $E^{k+1}_{i_2,\dots,i_\alpha,2p}\to E^{k+1}_{i_2,\dots,i_\alpha,2p+1}$ qui rend commutatif la totalité du diagramme et telle que le carré ainsi complété est cartésien.

Rappelons que l'espace $E^{k}_{i_2,\dots,i_\alpha}$ est défini comme la colimite homotopique de :
$$ \scriptstyle E^{k+1}_{i_2,\dots,i_\alpha,1}\gets E^{k+1}_{i_2,\dots,i_\alpha,2}\to E^{k+1}_{i_2,\dots,i_\alpha,3}
 \gets\cdots\to
 E^{k+1}_{i_2,\dots,i_\alpha,2\alpha-1}\gets E^{k+1}_{i_2,\dots,i_\alpha,2\alpha}\to E^{k+1}_{i_2,\dots,i_\alpha,2\alpha+1}.$$
Le diagramme commutatif précédent étant constitué de carrés cartésiens, le théorème \ref{Puppe} s'applique et le carré de la proposition :
$$\xymatrix{
E^{k+1}_{i_2,\dots,i_\alpha}\ar[rr]^{\phi^{k+1}_{i_2,\dots,i_\alpha}}\ar[d]	&& {\FF_{k+1} (A\times J)}\ar[d] \\
E^{k}_{i_2,\dots,i_\alpha}\ar[rr]^{\phi^{k}_{i_2,\dots,i_\alpha}}	&& {\FF_k (A\times J)} \\
}
$$
est cartésien, ce qui montre $H_s$.
\qed

Cette proposition nous permet de définir la notion de bord pour les espaces $E^{*}_{\ibar}$.
\Def Pour tout $\ibar=(i_2,\dots,i_\alpha)\in\Cal A_\alpha$ et tout entier $k$, $k\ge\alpha+1$, on définit ${\partial E^{k+1}_{\ibar}}$ comme le produit fibré suivant :
$$\xymatrix{
{\partial E^{k+1}_{\ibar}}\ar[rr]\ar[d]	&& {\partial\FF_{k+1} (A\times J)}\ar[d] \\
E^{k}_{\ibar}\ar[rr]^{\phi^{k}_{\ibar}}	&& {\FF_k (A\times J)} \\
}
$$

\Prop 
L'application $\phi^{k+1}_{\ibar}$ envoie la paire $({E^{k+1}_{\ibar}},{\partial E^{k+1}_{\ibar}})$ sur la paire $({\FF_{k+1}(A\times J)}, {\partial\FF_{k+1} (A\times J)})$. En particulier, le carré suivant est cartésien
$$\xymatrix{
{\partial E^{k+1}_{\ibar}}\ar[rr]^<(0.15){\phi^{k+1}_{\ibar}}\ar[d]	&& {\partial\FF_{k+1} (A\times J)=\FF_k (A\times J)}\times A\times \{-k,k\}\ar[d] \\
E^{k}_{\ibar}\ar[rr]_{\phi^{k}_{\ibar}}	&& {\FF_k (A\times J)} \\
}
$$
De plus, ${\partial E^{k+1}_{\ibar}}={E^{k}_{\ibar}}\times A\times \{-k,k\}$ et l'application 
${\partial E^{k+1}_{\ibar}} \stackrel{\phi^{k+1}_{\ibar}}{\longrightarrow} {\partial\FF_{k+1} (A\times J)}$
 est donnée par $\phi^{k}_{\ibar}\times Id_{A\times\{-k,k\}}$

\Dem La preuve découle immédiatement du fait que le produit fibré homotopique d'une fibration triviale est lui aussi trivial.\qed

Fixons un entier $k\ge 1$. En appliquant le résultat de la proposition \ref{pull} pour $\alpha=1$, on obtient un espace $E^k$ qui est la colimite de : 
$$E^k_1\gets E^k_2 \to E^k_3 .$$
Cet espace $E^k$ a le type d'homotopie de l'espace de configurations $\FF_k(A\times J)$. Plus précisément, on a :
\Th \label{ThConstruct} Pour tout $k\ge 1$, il existe un carré cartésien
$$\xymatrix{
E^{k+1}\ar[rr]^{\phi^{k+1}}\ar[d]	&& {\FF_{k+1} (A\times J)}\ar[d] \\
E^{k}\ar[rr]^{\phi^{k}}	&& {\FF_k (A\times J)} \\
}
$$
dans lequel les flèches ${\phi^{k+1}}$, ${\phi^{k}}$ sont une équivalence d'homotopie fibrée.

\Dem D'après la proposition \ref{pull}, le carré est cartésien pour tout $k$. Ainsi, si ${\phi^{k}}$ est une équivalence d'homotopie, il en est de même pour ${\phi^{k+1}}$. Rappelons que $\phi^1$ est une équivalence d'homotopie, une récurrence évidente termine alors la preuve.\qed

\section{Un invariant du type d'homotopie : $\FF_k (A\times\RR)$}
\label{sec:Invariant}

Dans ce paragraphe, nous allons démontrer le théorème A énoncé dans l'introduction. On se donne donc une $\Cal F_2$-équivalence d'homotopie $f:A\to B$ entre deux variétés.
Remarquons que cette hypothèse permet d'obtenir une équivalence d'homotopie entre $E^2(A)$ et $E^2(B)$ et donc entre $\FF_2(A\times\RR)$ et $\FF_2(B\times\RR)$ (cf. \ref{F2AxR}). Nous commençons cette section en exhibant deux exemples de conditions suffisantes pour avoir une $\Cal F_2$-équivalence d'homotopie entre deux variétés.

\Exp \label{F2equiv1} Considérons deux variétés $\check A$ et $\check B$ ayant même type d'homotopie et même dimension. Alors il existe une $\Cal F_2$-équivalence d'homotopie entre $A=\check A \times \RR$ et $B=\check B \times \RR$.

\Dem 
Soit $n$ la dimension commune à $\check A$ et $\check B$. Nous munissons $\check A$ (resp. $\check B$) d'une métrique riemannienne telle que pour tout point $a$ de $\check A$ la restriction de l'exponentielle au disque unitaire tangent en $a$, $\exp : \Cal D_a(\check A)\to\check A$, est injective (Ceci est possible d'après \cite{Greene1978}). Or, le fibré tangent de $A=\check A \times\RR$, noté $T_A$, est le produit du fibré tangent de $\check A$, noté $T_{\check A}$ et du fibré trivial au dessus de $\RR$, noté $\epsilon$. Notons $\sigma :\RR\to\epsilon$ la section unitaire du fibré $\epsilon$ telle que pour tout $t\in\RR$, $\exp(\sigma(t))=t+1$.  On note encore $\sigma:\check A\to T_A$ la section du fibré $T_{\check A}\times\epsilon$ obtenue à partir de la section nulle sur $\check A$ et $\sigma :\RR\to\epsilon$. Remarquons que $\sigma:\check A\to T_A$ est encore unitaire.

On dénote par $\Cal S_A$ le fibré tangent en sphère au dessus de $A$ et on note $\Cal S^+(\check A)$ le produit fibré suivant :
$$\xymatrix{
{\Cal S^+(\check A)}\ar[rr]\ar[d]	&&	{\Cal S _A}\ar[d]	\\
{\check A=\check A\times \{0\}}\ar[rr]\ar@/^/[u]^\sigma &&	{A=\check A\times\RR}\ar@/^/[u]^\sigma	\\
}$$

Considérons le diagramme suivant :
$$\xymatrix{
	& {\Cal S^+ (\check A)}\ar[rd]\ar[rrd]	\\
{\check A}\ar[ru]^\sigma\ar[rd]_\Delta	&	&	C_A\ar@{-->}[r]	&	{\FF_2 (A) }	\\
	& {\check A\times\check A}\ar[ru]\ar[rru]	\\
}$$
où $\Delta$ est l'application diagonale et le carré intérieur est une somme amalgamée.
L'application $ \Cal S^+ (\check A) \to	{\FF_2 (A) } $ est donnée par $ (a,v)\mapsto (a,\exp v)$ et ${\check A\times\check A}\to	{\FF_2 (A) }$ par $(a_1,a_2)\mapsto((a_1,0),(a_2,1))$.
Ces applications sont bien définies grâce aux choix effectués pour $\sigma$. On vérifie alors que le carré extérieur commute. Par propriété universelle d'une somme amalgamée, on obtient une application induite $C_A \to \FF_2 (A)$ qui fait commuter le diagramme. Nous appliquons alors \ref{Puppe2} avec comme base $\check A$ et montrons que $C_A \to \FF_2 (A)$ est une équivalence d'homotopie. Nous avons évidemment la même construction pour la variété $B$ ; nous allons donc exhiber une équivalence d'homotopie entre $C_A$ et $C_B$.

 Les résultats de R.~Benlian et J.~Wagoner \cite{Benlian1967} , montrent que si $\check A$ et $\check B$ sont deux variétés ayant même type d'homotopie et même dimension $n$, alors leurs fibrés tangents en sphère, ${\Cal S _{\check A}}$  et ${\Cal S _{\check B}}$, ont même type d'homotopie fibrée. Le diagramme suivant commute donc à homotopie près :
$$\xymatrix{
{\Cal S _{\check A}}\ar[rr]^\he\ar[d]		&&	{\Cal S _{\check B}}\ar[d]	\\
{\check A}\ar[rr]^\he		&&	{\check B}\\
}$$
Remarquons que $S^+(\check A)$ (resp. $S^+(\check B)$) est aussi la suspension suivant la fibre\footnote{«fiberwise suspension»} de $S_{\check A}$ (resp. $S_{\check B}$). Nous avons donc une équivalence d'homotopie
fibrée entre $\Cal S^+ (\check A)$ et $\Cal S^+ (\check B)$ qui, par la définition de $\sigma$, est
compatible avec les sections.
Plus précisément, le diagramme suivant commute à homotopie près
$$\xymatrix{
{\Cal S^+ (\check A)}\ar[rr]^\he\ar[d]		&&	{\Cal S^+ (\check B)}\ar[d]	\\
{\check A}\ar[rr]^\he\ar@/_/[u]_\sigma\ar[d]		&&	{\check B}\ar@/_/[u]_\sigma\ar[d]\\
{\check A\times\check A}\ar[rr]^\he		&&	{\check B\times\check B}\\
}$$
et les flèches horizontales sont des équivalences d'homotopie. On obtient alors une équivalence d'homotopie entre $C_A$ et $C_B$ qui définit une équivalence d'homotopie entre ${\FF_2 (A) }$ et ${\FF_2 (B) }$.

La commutativité à homotopie près du carré
$$ \xymatrix{
 {\FF_2 (A)}\ar[d]_\he\ar[r]	&	A\times A\ar[d]_{f\times f}^\he	\\
 {\FF_2 (B)}\ar[r]		&	B\times B	\\
 }
$$
s'obtient facilement en remarquant que  $\Cal S^+ (\check A) \to \FF_2 (A)$ se factorise à homotopie près en $\Cal S^+ (\check A) \to  \check A \to \check A\times \check A \to \FF_2 (A)$.\qed

\Rem Une construction utilisant la même technique que dans la démonstration précédente permet d'obtenir, sous les mêmes hypothèses, une équivalence d'homotopie entre $\FF_3 (\check A \times \RR)$ et $\FF_3 (\check B \times \RR)$.
 
\Exp \label{F2equiv2} Considérons deux variétés compactes et $2$-connexes, A et B, ayant le même type d'homotopie. Alors il existe une $\Cal F_2$-équivalence d'homotopie entre $A$ et $B$.

\Dem Le résultat découle de la démonstration du théorème A de \cite{Aouina2003}. Donnons une idée de la preuve dans le cas $3$-connexe en utilisant des résultats antérieurs de J.~R.~Klein \cite{Klein1999,Klein2002}. Soit $f:A\to B$ une équivalence d'homotopie. On note respectivement ${\Cal S_A}$, ${\Cal S_B}$ les fibrés tangents en sphère au dessus de $A$ et $B$ et $n$ la dimension commune à $A$ et $B$. On considère maintenant le diagramme suivant :
$$\xymatrix{
{\Cal S_A}\ar[rrr]^\exp\ar[ddd]	&			&		&	{\FF_2(A)}\ar[ddl]\ar[ddd]\\
			&{\Cal S_B}\ar[r]^\exp\ar[d]\ar[ul]_\he^\psi	&{\FF_2(B)}\ar[d]\ar@{.}[ru]	\\
			&B\ar[r]^\Delta\ar[ld]_{f^{-1}}			&B\times B	\\
A	\ar[rrr]	&	&	&	A\times A\ar[lu]^{f\times f}		\\
}$$
L'existence d'une équivalence d'homotopie $\psi:{\Cal S_A}\to{\Cal S_B}$ rendant le carré de gauche commutatif est garantie par les résultats de R.~Benlian et J.~Wagoner \cite{Benlian1967}.

On équipe $B$ (resp. $A$) d'une métrique telle que l'exponentielle soit définie et injective sur chaque espace tangent jusqu'au rayon $1$. La flèche $\exp:{\Cal S_B}\to\FF_2(B)$ où $\exp(b,v)=(b,\exp v)$ est alors bien définie.

Le carré de sommets $(B,\Cal S_B,\FF_2(B),B\times B)$ est cocartésien. En effet, le fibré normal de l'application diagonale $\Delta:B\to B\times B$ est le fibré tangent en sphère de $B$.
De même, le carré de sommets $(A,\Cal S_A,\FF_2(A),A\times A)$ est cocartésien et donc celui 
de sommets $(B,\Cal S_B,\FF_2(A),B\times B)$ l'est aussi.
On a ainsi obtenu deux complémentaires homotopiques à dualité de Poincaré de la diagonale $\Delta : B \to B\times B$. Nous appliquons alors le corollaire B de \cite{Klein2002} pour en déduire l'équivalence de $\FF_2(A)$ et $\FF_2(B)$ au dessus de $B\times B$. Pour cela, on remarque que :\\
\textbullet\quad $B$ a le type d'homotopie d'un CW-complexe de dimension $n$.\\
\textbullet\quad $B\times B$  a le type d'homotopie d'un CW-complexe de dimension $2n$.\\
\textbullet\quad $\Delta : B \to B\times B$ est $3$-connexe avec $3 > \mathrm{dim}(B\times B) - 2\mathrm{dim}(B) + 2$.
\qed

\Prop \label{equiv1} Considérons deux variétés $A$ et $B$ reliées par une $\Cal F_2$-équivalence d'homotopie. Pour tout $k\ge 1$ et tout $\ibar=(i_2,\dots,i_k)\in\Cal A_k$, il existe une équivalence d'homotopie fibrée :
$$\xymatrix{
E^{k+1}_{\ibar}(A)\ar[r]^\he_{\Psi^{k+1}_{\ibar}}\ar[d]	&E^{k+1}_{\ibar}(B)\ar[d]\\
E^{k}_{\ibar}(A)\ar[r]_\he^{f^{\times k}}		&E^{k}_{\ibar}(B)\\
}$$
où $\Psi^{k+1}_{\ibar}$ se restreint à l'application 
$$f^{\times k}\times f\times Id : \partial E^{k+1}_{\ibar}(A)=E^{k}_{\ibar}(A)\times A\times\{-k,k\}\stackrel{\he}{\longrightarrow} \partial E^{k+1}_{\ibar}(B)=E^{k}_{\ibar}(B)\times B\times \{-k,k\} .$$

\Dem Soit $k\ge 1$ et $\ibar=(i_2,\dots,i_k)\in\Cal A_k$.

Nous utilisons le diagramme suivant dans lequel les flèches verticales sont des équivalences d'homotopie :
$$\xymatrix@C=1.5ex{
E^{k+1}_{\ibar 1}(A)\ar[d] & E^{k+1}_{\ibar 2}(A)\ar[d]\ar[l]\ar[r] & E^{k+1}_{\ibar 3}(A)\ar[d]& \ar[l] {\dots}\ar[r] & E^{k+1}_{\ibar (2k-1)}(A)\ar[d]& E^{k+1}_{\ibar (2k)}(A)\ar[d]\ar[l]\ar[r] & E^{k+1}_{\ibar (2k+1)}(A)\ar[d]\\
E^{k+1}_{\ibar 1}(B) & E^{k+1}_{\ibar 2}(B)\ar[l]\ar[r] & E^{k+1}_{\ibar 3}(B) & \ar[l] {\dots}\ar[r] & E^{k+1}_{\ibar (2k-1)}(B)& E^{k+1}_{\ibar (2k)}(B)\ar[l]\ar[r] & E^{k+1}_{\ibar (2k+1)}(B)\\
}$$

Rappelons de la preuve de la proposition \ref{pull1} que chaque carré commute à homotopie près et équivaut au carré suivant :
$$\xymatrix{
A^{k+1}_{k*}\cong A^{k-1}\times\FF_2 (A)\ar[r]\ar[d]^\he	& A^{k+1}\cong A^{k-1}\times A^2\ar[d]^\he\\
B^{k+1}_{k*}\cong B^{k-1}\times\FF_2 (B)\ar[r]	& B^{k+1}\cong B^{k-1}\times B^2\\
}$$
pour un certain entier $*$, $1\le*<k$. On déduit donc une équivalence d'homotopie $\Psi$ entre la colimite de la première ligne et la colimite de la deuxième c'est-à-dire entre $E^{k+1}_{\ibar}(A)$ et $E^{k+1}_{\ibar}(B)$.

Rappelons que l'espace $\partial{E}^{k+1}_{\ibar}(A)$ a été défini comme le produit fibré suivant~:
$$\xymatrix@C=1.5ex{
 {\partial E^{k+1}_{\ibar}(A)}\ar[rr]^-{\phi^{k+1}_{\ibar}}\ar[d]	&&	{\partial\FF_{k+1} (A\times J)}=\FF_k (A\times J)\times (A\times  \{-k,k\})\ar[d]\\
 E^k_{\ibar}(A)= A^k\ar[rr]_{\phi^k_{\ibar}}			&&	{\FF_k (A\times J)}\\
}$$
D'après la définition de $\phi^{k+1}_{\ibar}$, on a :
$${\partial E^{k+1}_{\ibar}(A)}=E^{k+1}_{\ibar 1}(A)\coprod E^{k+1}_{\ibar (2k+1)}(A)=E^k_{\ibar}\times(A\times  \{-k,k\})$$
ce qui montre que l'application $\Psi$ se restreint bien à l'application 
$$f^{\times k}\times f\times Id : \partial E^{k+1}_{\ibar}(A)\stackrel{\he}{\longrightarrow} \partial E^{k+1}_{\ibar}(B) .$$

Finalement, pour tout $p\in[\![1,k]\!]$, nous avons la commutativité à homotopie près du diagramme~:
$$\xymatrix{
&E^{k+1}_{\ibar,2p}(B)\ar[rrr]\ar[dd]	&	&	&E^{k+1}_{\ibar,2p+ 1}(B)\ar[dd]\\
E^{k+1}_{\ibar,2p}(A)\ar[rrr]\ar[dd]\ar[ru]^\he	&	&	&E^{k+1}_{\ibar,2p+ 1}(A)\ar[dd]\ar[ru]^\he\\
&E^{k}_{\ibar}(B)=B^{\times k}\ar@{=}[rrr]	&	&	&E^{k}_{\ibar}(B)=B^{\times k}\\
E^{k}_{\ibar}(A)=A^{\times k}\ar@{=}[rrr]\ar[ru]_{f^{\times k}}^\he	&	&	&E^{k}_{\ibar}(A)=A^{\times k}\ar[ru]_{f^{\times k}}^\he\\
}$$
et d'un diagramme analogue avec $E^{k+1}_{\ibar,2p-1}(\bullet)$.
La commutativité du carré de la proposition s'ensuit.
\qed

 On a donc une équivalence d'homotopie $E^k_{\ibar}(A)\to E^k_{\ibar}(B) $ pour $\ibar\in\Cal A_k \bigcup\Cal A_{k+1}$. Comme dans la construction des $E^k_{\ibar}$ , on va effectuer une récurrence sur l'entier $k-\alpha$ afin d'obtenir des équivalences d'homotopie $E^k_{\ibar}(A)\to E^k_{\ibar}(B)$ pour $\ibar\in\Cal A_\alpha$ avec $1\le \alpha \le k$.

\Prop \label{equiv} Soit $\alpha \ge 1$ et $\ibar=(i_2,\dots,i_\alpha)\in\Cal A_\alpha$ alors, pour tout entier $k$, $k\ge\alpha$, il existe des équivalences d'homotopie  ${\Psi^{k}_{\ibar}}:(E^{k}_{\ibar}(A),\partial E^{k}_{\ibar}(A))\to (E^{k}_{\ibar}(B),\partial E^{k}_{\ibar}(B))$  telles que l'on ait l'équivalence d'homotopie fibrée suivante :
$$\xymatrix{
E^{k+1}_{\ibar}(A)\ar[rr]_\he^{\Psi^{k+1}_{\ibar}}\ar[d]	&& E^{k+1}_{\ibar}(B)\ar[d] \\
E^{k}_{\ibar}(A)\ar[rr]^\he_{\Psi^{k}_{\ibar}}	&& E^{k}_{\ibar}(B) \\}
$$
et la restriction de l'application $\Psi^{k+1}_{\ibar}$ à l'espace $\partial E^{k+1}_{\ibar}(A)$ et à valeurs dans $\partial E^{k+1}_{\ibar}(B)$ est homotope à l'application
$$\Psi^{k}_{\ibar}\times f\times Id : \partial E^{k+1}_{\ibar}(A)=E^{k}_{\ibar}(A)\times A\times\{-k,k\}\stackrel{\he}{\longrightarrow}\partial E^{k+1}_{\ibar}(B)=E^{k}_{\ibar}(B)\times B\times\{-k,k\} .$$

\Dem Comme dans la preuve de la proposition \ref{pull}, nous procédons par récurrence sur l'entier $k-\alpha$.
Pour tout entier $s$, $s\ge 1$, nous définissons la propriété $H_s$ par :
pour tout couple d'entiers $(k,\alpha)$, $1\le k-\alpha\le s$, et tout $\ibar\in\Cal A_\alpha$, il existe des équivalences d'homotopie
$$\Psi^{k+1}_{\ibar,*} : (E^{k+1}_{\ibar,*}(A),\partial E^{k+1}_{\ibar,*}(A))\to (E^{k+1}_{\ibar,*}(B),\partial E^{k+1}_{\ibar,*}(B)) \textrm{ où } *\in \{1,3,\dots,2\alpha+1\}$$
$$\textrm{ et } \Psi^{k}_{\ibar,*} : E^{k}_{\ibar,*}(A)\to E^{k}_{\ibar,*}(B) \textrm{ où } *\in [\![1,2\alpha+1]\!] $$
telles que pour tout $p\in[\![1,\alpha]\!]$, le diagramme  suivant, noté $\Cal D^k_{\ibar}(p)$, commute à homotopie près :
$$\xymatrix@C=1.5ex{
& E^{k+1}_{\ibar (2p-1)}(B)\ar[dd] && E^{k+1}_{\ibar (2p)}(B)\ar[ll]\ar[rr]\ar[dd] && E^{k+1}_{\ibar (2p+1)}(B)\ar[dd]\\
E^{k+1}_{\ibar (2p-1)}(A)\ar[dd]\ar[ur]^\he && E^{k+1}_{\ibar (2p)}(A)\ar[ll]\ar[rr]\ar[dd] && E^{k+1}_{\ibar (2p+1)}(A)\ar[dd]\ar[ur]^\he\\
&E^{k}_{\ibar (2p-1)}(B)&& E^{k}_{\ibar (2p)}(B)\ar[ll]\ar[rr] && E^{k}_{\ibar (2p+1)}(B)\\
E^{k}_{\ibar (2p-1)}(A)\ar[ur]_\he&& E^{k}_{\ibar (2p)}(A)\ar[ll]\ar[rr]\ar[ur]_\he && E^{k}_{\ibar (2p+1)}(A)\ar[ur]_\he\\
}$$
et telles que la restriction de l'application $\Psi^{k+1}_{\ibar,*}$ à l'espace $\partial E^{k+1}_{\ibar,*}(A)$ et à valeurs dans $\partial E^{k+1}_{\ibar,*}(B)$ est homotope à l'application
$$\Psi^{k}_{\ibar,*}\times f\times Id : \partial E^{k+1}_{\ibar,*}(A)=E^{k}_{\ibar,*}(A)\times A\times\{-k,k\}\stackrel{\he}{\longrightarrow}\partial E^{k+1}_{\ibar,*}(B)=E^{k}_{\ibar,*}(B)\times B\times\{-k,k\}.$$

La propriété $H_1$ découle de la proposition \ref{equiv1}. En effet, si $k=\alpha+1$ et $\ibar\in\Cal A_{\alpha}$, les carrés droit et gauche du diagramme sont la conclusion de la proposition \ref{equiv1} et les carrés inférieurs ont été explicités dans la démonstration.

Supposons dès lors que pour un certain entier $s$, $s\ge 2$, la propriété $H_{s-1}$ est vraie.
On fixe $k=\alpha+s$, $\ibar=(i_2,\dots,i_{\alpha})\in\Cal A_{\alpha}$ et $p\in [\![1,\alpha]\!]$.
D'après l'hypothèse de récurrence $H_{s-1}$, on a des équivalences d'homotopie
$E^{k}_{\ibar,*}(A)\to E^{k}_{\ibar,*}(B)$ pour $*\in \{1,3,\dots,2\alpha+1\}$.

Par définition des espaces $E^*_{\ibar (2p)}(\bullet)$, les faces avant et arrière du diagramme $\Cal D^{k-1}_{\ibar}(p)$ sont constituées de carrés cartésiens. Le lemme \ref{bicube} permet alors d'obtenir l'équivalence d'homotopie 
$$\Psi^{k}_{\ibar,(2p)} : E^{k}_{\ibar,(2p)}(A)\to E^{k}_{\ibar,(2p)}(B)$$ qui commute avec le reste du diagramme $\Cal D^{k-1}_{\ibar}(p)$.
En remarquant que la partie supérieure du diagramme $\Cal D^{k-1}_{\ibar}(p)$ est exactement la partie inférieure de $\Cal D^{k}_{\ibar}(p)$, nous obtenons la commutativité de la partie inférieure de $\Cal D^{k}_{\ibar}(p)$.
Ce même argument appliqué aux diagrammes $\Cal D^{k}_{\ibar,2p+1}(q)$, $q\in [\![1,\alpha+1]\!]$, permet d'obtenir l'equivalence d'homotopie 
$$\Psi^{k}_{\ibar,(2p+1),2q} : E^{k}_{\ibar,(2p+1),2q}(A)\to E^{k}_{\ibar,(2p+1),2q}(B)$$
qui commute à homotopie près avec les applications du diagramme $\Cal D^{k}_{\ibar,2p+1}(q)$.

Rappelons que l'espace $E^{k}_{\ibar,(2p+1)}(A)$ est défini comme la colimite homotopique de :
$$\scriptstyle E^{k}_{\ibar,(2p+1),1}(A)\gets E^{k}_{\ibar,(2p+1),2}(A)\to E^{k}_{\ibar,(2p+1),3}(A)
 \gets\cdots\to
 E^{k}_{\ibar,(2p+1),2\alpha+1}(A)\gets E^{k}_{\ibar,(2p+1),2\alpha+2}(A)\to E^{k}_{\ibar,(2p+1),2\alpha+3}(A)$$
et que les espaces $E^{k}_{\ibar,(2p+1)}(B)$,  $E^{k+1}_{\ibar,(2p+1)}(A)$ et  $E^{k+1}_{\ibar,(2p+1)}(B)$ sont définis de façon similaire. 
 
Les diagrammes $D^{k}_{\ibar,(2p+1)}(q)$, $q\in [\![1,\alpha+1]\!]$, complétés par les applications $\Psi^{k}_{\ibar,(2p+1),2q}$ induisent donc, en passant à la colimite selon $q$, une équivalence d'homotopie  $\Psi^{k+1}_{\ibar,(2p+1)}$ qui fait commuter le  diagramme suivant :
$$\xymatrix{
E^{k+1}_{\ibar,(2p+1)}(A)\ar[rr]^\he_{\Psi^{k+1}_{\ibar,(2p+1)}}\ar[d]	&& E^{k+1}_{\ibar,(2p+1)}(B)\ar[d] \\
E^{k}_{\ibar,(2p+1)}(A)\ar[rr]_\he^{\Psi^{k}_{\ibar,(2p+1)}}	&& E^{k}_{\ibar,(2p+1)}(B) \\}
$$
Nous avons donc complété le carré droit du diagramme $\Cal D^{k}_{\ibar}(p)$.
On complète le carré gauche de manière analogue.

En remplaçant, dans la partie supérieure du diagramme $\Cal D^{k}_{\ibar}(p)$, les espaces $E^{k+1}_{\ibar,*}(\bullet)$ par leurs bords $\partial E^{k+1}_{\ibar,*}(\bullet)$ et en utilisant la propriété \ref{bicube}, nous obtenons une application
$$^\partial \Psi^{k+1}_{\ibar,(2p)} :\partial E^{k+1}_{\ibar,(2p)}(A)\to\partial E^{k+1}_{\ibar,(2p)}(B) .$$
Rappelons de la démonstration de la propriété \ref{bicube}, que les applications $^\partial \Psi^{k+1}_{\ibar,(2p)}$ et $\Psi^{k+1}_{\ibar,(2p)}$ sont obtenues par la propriété universelle du produit fibré homotopique le long de l'application $E^{k}_{\ibar,(2p+1)}(A)\to E^{k}_{\ibar,(2p+1)}(B)$ qui apparaît dans le carré droit du diagramme commutatif, à homotopie près, suivant :
$$\xymatrix{
\partial E^{k+1}_{\ibar,(2p+1)}(A)\ar[r]\ar[d]&	E^{k+1}_{\ibar,(2p+1)}(A)\ar[d]_{\Psi^{k+1}_{\ibar,(2p+1)}}\ar[r]	&  E^{k}_{\ibar,(2p+1)}(A)\ar[d]_{\Psi^{k}_{\ibar,(2p+1)}}	&	E^{k}_{\ibar,(2p)}(A)\ar[l]\ar[d]^{\Psi^{k}_{\ibar,(2p)}}\\
\partial E^{k+1}_{\ibar,(2p+1)}(B)\ar[r]	&	E^{k+1}_{\ibar,(2p+1)}(B)\ar[r]	& 
 E^{k}_{\ibar,(2p+1)}(B)	&	E^{k}_{\ibar,(2p)}(B)\ar[l]\\
}$$
Cela nous assure que la restriction de l'application $\Psi^{k+1}_{\ibar}$ à l'espace $\partial E^{k+1}_{\ibar}(A)$, et à valeurs dans $\partial E^{k+1}_{\ibar}(B)$, est homotope à l'application $^\partial \Psi^{k+1}_{\ibar,(2p)}$.

Rappelons aussi qu'un produit fibré homotopique obtenu à partir d'une fibration triviale est lui aussi trivial. L'application $^\partial \Psi^{k+1}_{\ibar,(2p)}$ est donc égale à l'application 
$${\Psi^{k}_{\ibar,(2p)}}\times f\times Id : \partial E^{k+1}_{\ibar,(2p)}(A)\stackrel{\he}{\longrightarrow} \partial E^{k+1}_{\ibar,(2p)}(B)$$
Nous en déduisons, d'après \ref{trivcolim}, que la restriction de l'application $\Psi^{k+1}_{\ibar}$ à l'espace $\partial E^{k+1}_{\ibar}(A)$, et à valeurs dans $\partial E^{k+1}_{\ibar}(B)$, est homotope à l'application
$$\Psi^{k}_{\ibar}\times f\times Id : \partial E^{k+1}_{\ibar}(A)=E^{k}_{\ibar}(A)\times A\times\{-k,k\}\stackrel{\he}{\longrightarrow}\partial E^{k+1}_{\ibar}(B)=E^{k}_{\ibar}(B)\times B\times\{-k,k\} .$$
D'où la propriété $H_{s}$.
\qed

Le résultat suivant contient le théorème A de l'introduction : 
\Th \label{EquivBords}
Soit $f:A\to B$ une $\Cal F_2$-équivalence d'homotopie entre deux variétés $A$ et $B$, alors, pour tout $k\ge 1$, on a une équivalence d'homotopie fibrée :
$$\xymatrix{
({\FF_{k+1} }(A\times J),\partial{\FF_{k+1} }(A\times J))\ar[rr]_\he^{f_{k+1}}\ar[d]	&& ({\FF_{k+1} }(B\times J),\partial{\FF_{k+1} }(B\times J))\ar[d] \\
{{\FF_k (A\times J)}}\ar[rr]_\he^{f_k}	&& {{\FF_k (B\times J)}} \\}
$$
où $f_{k+1} : \partial{\FF_{k+1} }(A\times J)\to \partial{\FF_{k+1} }(B\times J)$ est donnée, à homotopie près, par le produit cartésien suivant :
$$f_k \times (f\times Id) : {\FF_k (A\times J)}\times (A\times \{-k,k\}) \to {\FF_k (B\times J)}\times (B\times \{-k,k\}) .$$

\Dem On utilise la proposition \ref{equiv} dans le cas $\alpha=1$ pour obtenir les équivalences d'homotopie $f_{k+1}$ et $f_{k}$ :
$$\xymatrix{
{\FF_{k+1} }(A\times  J)\ar[r]^\he\ar[d]	&E^{k+1}(A)\ar[r]^\he\ar[d]	& E^{k+1}(B)\ar[r]^\he\ar[d]& {\FF_{k+1} }(B\times  J)\ar[d] \\
{\FF_k }(A\times  J)\ar[r]_\he	&E^{k}(A)\ar[r]_\he & E^{k}(B)\ar[r]_\he& {\FF_k }(B\times  J) \\}
$$
L'assertion sur $f_{k+1} : \partial{\FF_{k+1} }(A\times J)\to \partial{\FF_{k+1} }(B\times J)$ est évidente en utilisant la définition de $\partial{\FF_{k+1} }(\bullet\times J)$ et le résultat précédent.
\qed

\Cor \label{PullBord}
On considère deux espaces $P_A$ et $P_B$. Si $\mu:P_A\to P_B$ est une équivalence d'homotopie telle que le carré
$$\xymatrix{
P_A\ar[rr]\ar[d]^\he_\mu	&&	{\FF_k (A\times J)}\ar[d]_\he^{f_k}	\\
P_B\ar[rr]	&&	{\FF_k (B\times J)}	\\
}$$
commute à homotopie près, alors on a un carré commutatif à homotopie près :
$$\xymatrix{
P_A \times (A\times\{-k,k\})\ar[rr]\ar[d]^\he_{\mu\times f\times Id}	&&	{\FF_{k+1} (A\times J)}\ar[d]_\he^{f_{k+1}}\\
P_B \times (B\times\{-k,k\})\ar[rr]	&&	{\FF_{k+1} (B\times J)}	\\
}$$

\Dem 
La commutativité à homotopie près du carré de l'hypothèse nous permet d'obtenir la commutativité à homotopie près du carré :
$$\xymatrix@C=1.5ex{
P_A \times (A\times\{-k,k\})\ar[rr]\ar[d]^\he_{\mu\times f}&&	{\FF_{k} (A\times J)\times(A\times \{-k,k\})}\ar[d]^\he\\
P_B \times (B\times\{-k,k\})\ar[rr]&&	{\FF_{k} (B\times J)\times(B\times \{-k,k\})}\\
}$$

Or ${\FF_{k} (A\times J)\times(A\times \{-k,k\})}={\partial\FF_{k+1} (A\times J)}$, on a donc montré que le carré de gauche du diagramme suivant est commutatif à homotopie près :
$$\xymatrix@C=1.5ex{
P_A \times (A\times\{-k,k\})\ar[rr]\ar[d]^\he_{\mu\times f}&&	{\partial\FF_{k+1} (A\times J)}\ar[d]^\he\ar[rr]	&&	{\FF_{k+1} (A\times J)}\ar[d]^\he	\\
P_B \times (B\times\{-k,k\})\ar[rr]&&	{\partial\FF_{k+1} (B\times J)}\ar[rr]	&&	{\FF_{k+1} (B\times J)}\\
}$$
La proposition précédente nous assure que le carré de droite de ce diagramme commute à homotopie près, d'où le corollaire.\qed

\section{Un invariant du type d'homotopie : $\Sigma^{k-2}\ \FF_k (A)$}
\label{sec:Suspension}

Dans ce paragraphe, nous allons montrer que si deux variétés connexes $A$ et $B$ ont même $\Cal F_2$-type d'homotopie  alors une certaine suspension des espaces $\FF_2(A)$ et $\FF_2(B)$ ont même type d'homotopie. On supposera toujours que la dimension des variétés est au moins égale à trois, les résultats restant d'ailleurs vrais en dimension inférieure.

\Th \label{ThSusp} Si l'on a une $\Cal F_2$-équivalence d'homotopie $f:A\to B$  entre deux variétés connexes $A$ et $B$, alors les suspensions $\Sigma^{k-2} \FF_k (A)$ et $\Sigma^{k-2}\FF_k (B)$ ont même type d'homotopie.

La preuve du théorème \ref{ThSusp} utilise la caractérisation du terme $E^k_2(A)$ comme produit fibré homotopique :
$$\xymatrix{
E_2^k(A)	\ar[r]\ar[d]	&	{\FF_k (A\times J)}\ar[d]\\
E_2^{k-1}(A)\ar[r]	&	{\FF_{k-1} (A\times J)}\\
}$$

\Prop \label{ThSusp1} Soit $A$ et $B$ deux variétés connexes ayant même $\Cal F_2$-type d'homotopie. Alors les suspensions $\Sigma \FF_3 (A)$ et $\Sigma\FF_3 (B)$ ont même type d'homotopie.

\Dem 
Notons $C$ la colimite homotopique de : 
$$\FF_3 (A) \gets{\FF_3 (A)\times\partial D^1}\to {\FF_2 (A)\times (A\times\partial D^1)}$$
où $D^1$ est le disque unité de $\RR$. En projetant chaque espace sur $\FF_2(A)$ et en appliquant \ref{Puppe}, on constate que l'application $C\to \FF_2(A)$ a pour fibre homotopique $A\times\{-1,1\}\bigcup (A\setminus\{q_1,q_2\})\times [-1,1]$.

Considérons les applications 
$$\begin{array}{rcl}
{\FF_3 (A)} &\to&{\FF_3  (A\times J)} \\
(a_1,a_2,a_3)&\mapsto& ((a_1,0),(a_2,0),(a_3,0))\\
\end{array}
$$ et
$$\begin{array}{rcl}
{\FF_2 (A)\times (A\times\partial D^1)}&\to&{\FF_3  (A\times J)} \\
(a_1,a_2,a_3,-1)&\mapsto& ((a_1,0),(a_2,0),(a_3,-2))\\
(a_1,a_2,a_3,+1)&\mapsto& ((a_1,0),(a_2,0),(a_3,+2))\\
\end{array}$$
qui font commuter les carrés suivants :
$$\xymatrix{
{\FF_3 (A)}\ar[r]\ar[d]&{\FF_3  (A\times J)}\ar[d]	&&{\FF_2 (A)\times (A\times\partial D^1)}\ar[r]\ar[d]&{\FF_3  (A\times J)}\ar[d]	\\
{\FF_2 (A)}\ar[r]&{\FF_2  (A\times J)}	&&{\FF_2 (A)}\ar[r]					&{\FF_2  (A\times J)}	\\
}$$
Par la propriété universelle d'un produit fibré homotopique, ces deux applications factorisent au travers de l'application $E^3_2(A)\to{\FF_3  (A\times J)}$. On en déduit alors la commutativité à homotopie près du carré extérieur du diagramme suivant :
$$ \xymatrix@C=1.5ex{
          	& {\FF_3 (A)}\ar@/^/[rrrrd]\ar[rrd]             			\\
{\FF_3 (A)\times\partial D^1}\ar[ur]\ar[rrd]	&         	&         	& C\ar@{-->}[rr]	& &{E^3_2(A)}	\\
        	&	&{\FF_2 (A)\times (A\times\partial D^1)}\ar@/_/[urrr]\ar[ur]	&	\\		
}
$$
Avec la propriété universelle de la somme amalgamée homotopique, on en déduit un diagramme commutatif à homotopie près :
$$\xymatrix{
C\ar[r]\ar[d]		&	{E^3_2(A)}\ar[d]\ar[r]	&{\FF_3  (A\times J)}\ar[d]\\
{\FF_2 (A)}	\ar@{=}[r]	&	{E^2_2(A)}	\ar[r]	&{\FF_2  (A\times J)}\\
}$$
D'aprés \ref{pull}, le carré de droite est cartésien. De plus, la restriction de l'application composée $C \to {E^3_2(A)}\to {\FF_3  (A\times J)}$ aux fibres homotopiques des deux flèches verticales du grand carré est donnée par l'application~:
$$(A\times\{-2,2\})\bigcup (A\setminus\{q_1,q_2\})\times [-2,2]\hookrightarrow (A\times[-2,2])\setminus\{(q_1,0),(q_2,0)\}$$
qui est une équivalence d'homotopie. D'aprés \cite{May1975},  le grand carré du diagramme précédent est cartésien.
Le lemme du prisme (\ref{lemprism}) permet de conclure que $C$ a le même type d'homotopie que $E^3_2(A)$. En particulier, la cofibre de ${\FF_2 (A)\times (A\times\partial D^1)} \to E_2^3(A)$ a même type d'homotopie que la cofibre de ${\FF_3 (A)\times\partial D^1}\to\FF_3 (A)$ qui est $\Sigma \FF_3 (A)_+$ : la suspension de l'union disjointe de $\FF_3(A)$ avec un point.

Remarquons aussi que le terme $ {\FF_2 (A)\times (A\times\partial D^1)} $ correspond à l'espace $\partial E_2^3$. Or, d'après \ref{PullBord} le diagramme suivant commute :
$$\xymatrix{
{\FF_2 (A)\times (A\times\partial D^1)}\ar[r]\ar[d]_\he	&	E_2^3(A)\ar[d]_\he\ar[r] & {\FF_3  (A\times J)}\ar[d]_\he\\
{\FF_2 (B)\times (B\times\partial D^1)}\ar[r]	& E_2^3(B)\ar[r] & {\FF_3  (B\times J)}\\
}$$
On en déduit donc une équivalence d'homotopie entre $\Sigma \FF_3 (A)_+$ et $\Sigma \FF_3 (B)_+$. 
Or $\FF_3 (A)$ et $\FF_3 (B)$ sont connexes car $A$ et $B$ sont connexes et de dimension supérieure à trois.
L'existence de cette équivalence d'homotopie implique celle d'une équivalence d'homotopie entre $\Sigma \FF_3 (A)$ et $\Sigma \FF_3 (B)$. 
On pourra trouver la démonstration de cette affirmation dans \cite[lemme 2.2]{Aouina2003}\qed

\Demr \ref{ThSusp}:

La démonstration s'inspire du cas $k=3$. Nous allons définir des espaces $ Z^k(A)$ tels que $E_2^k(A)$ soit la colimite homotopique de 
$$\FF_k (A) \gets{\FF_k (A)\times\partial D^{k-2}}\to  Z^k(A)$$
où $D^{k-2}$ est le disque unité de $\RR^{k-2}$. (On fait de même pour $B$).

Posons $ Z^3(A)={\FF_2 (A)\times (A\times\partial D^1)}$. Lors de la démonstration de la proposition \ref{ThSusp1}, nous avons vu que la colimite de 
$$\FF_3 (A) \gets{\FF_3 (A)\times\partial D^{1}}\to  Z^3(A)$$
est bien $E_2^3(A)$.

Nous allons maintenant prouver notre résultat par récurrence sur l'entier $k$.
On se donne donc un entier $k$ tel que :

\begin{description}
	\item [P1] : $E_2^k(A)$ (resp. $E_2^k(B)$) est la colimite homotopique de 
$$\FF_k (A) \gets{\FF_k (A)\times\partial D^{k-2}}\to Z^k(A) .$$
$$(\textrm{resp. } \FF_k (B) \gets{\FF_k (B)\times\partial D^{k-2}}\to Z^k(B)  .)$$
	\item [P2] : Il existe une équivalence d'homotopie entre $Z^k(A)$ et $Z^k(B)$ telle que le diagramme suivant commute à homotopie près :
$$\xymatrix{
Z^k(A)\ar[r]\ar[d]_\he	&	E_2^k(A)\ar[d]^\he	\\
Z^k(B)\ar[r]		&	E_2^k(B)	\\
}$$
	\item [P3] : La composée $${\FF_k (A)\times\partial D^{k-2}}\to Z^k(A)\to E_2^k(A)\to{\FF_k (A\times J)}\to A^{\times k}$$
coïncide avec la composée de l'injection canonique et de la projection canonique  
 $${\FF_k (A)\times\partial D^{k-2}}\to\FF_k (A) \hookrightarrow A^{\times k} .$$ (De même pour $B$). 
\end{description}

Introduisons aussi la notation $ \pull X $. Pour tout espace $X$ et toute application $X\to\FF_k (M)$ (pour un certain $M$), $ \pull X $ dénote le produit fibré homotopique suivant~:
$$\xymatrix{
{\pull X }\ar[d]\ar[r]	&	{\FF_{k+1} (M)}\ar[d] \\
X\ar[r]				&	{\FF_k (M)}	\\
}$$
En particulier, en utilisant l'application composée $Z^k(A)\to E_2^k(A)\to{\FF_k (A\times J)}$, ${\pull {Z^k(A)} }$ est le produit fibré homotopique :
$$\xymatrix{
{\pull {Z^k(A)} }\ar[d]\ar[r]	&	{\FF_{k+1} (A\times J)}\ar[d] \\
Z^k(A)\ar[r]				&	{\FF_k (A\times J)}	\\
}$$
En appliquant le lemme du prisme (\ref{lemprism}) avec comme base le triangle commutatif suivant~:
$$\xymatrix{
Z^k(A)\ar[rr]\ar[rd]		&&	{\FF_k (A\times J)}	\\
		&E_2^k(A)\ar[ru]	\\
}$$
on voit que le carré ci-dessous est cartésien :
$$\xymatrix{
{\pull {Z^k(A)} }\ar[d]\ar[r]	&	{E_2^{k+1}(A)}\ar[d] \\
Z^k(A)\ar[r]				&	{E_2^k(A)}	\\
}$$
Dans le cube ci-dessous, le carré du bas provient de l'hypothèse de récurrence {\bf P2}, le carré de droite est le résultat de la proposition \ref{equiv} :
$$\xymatrix{
&{\pull {Z^k(B)} }\ar[dd]\ar[rr]	&&	{E_2^{k+1}(B)}\ar[dd] \\
{\pull {Z^k(A)} }\ar[dd]\ar[rr]\ar@{-->}[ru]^\he	&&	{E_2^{k+1}(A)}\ar[dd]\ar[ru]^\he \\
&Z^k(B)\ar[rr]				&&	{E_2^k(B)}	\\
Z^k(A)\ar[rr]\ar[ru]^\he				&&	{E_2^k(A)}\ar[ru]^\he	\\
}$$
Nous avons alors une équivalence d'homotopie entre ${\pull {Z^k(A)} }$ et ${\pull {Z^k(B)} }$ qui fait commuter la totalité du cube à homotopie près.

Définissons $ Z^{k+1}(A) $ comme la colimite homotopique de 
$$ \pull {Z^k(A)} \gets Z^k(A) \times (A\times\partial D^1)\to E_2^k(A)\times (A\times\partial D^1)$$
où $Z^k(A) \times (A\times\partial D^1)\to{\pull {Z^k(A)} }$ est obtenue par propriété universelle du produit fibré homotopique :
$$\xymatrix{
Z^k(A) \times (A\times\partial D^1)\ar@/^/ [rrd]\ar@{-->}[rd]\ar@/_/[ddr] \\
&{\pull {Z^k(A)} }\ar[d]\ar[r]	&	{E_2^{k+1}(A)}\ar[d] \\
&Z^k(A)\ar[r]				&	{E_2^k(A)}\\
}$$
La flèche du haut étant la composée :
$$Z^k(A) \times (A\times\partial D^1)\to E_2^k(A)\times (A\times\partial D^1) \to {E_2^{k}(A)\times (A\times\{-k,k\})} \to {E_2^{k+1}(A)}$$
nous avons une équivalence d'homotopie entre les deux suites :
$$ \pull {Z^k(A)} \gets Z^k(A) \times (A\times\partial D^1)\to E_2^k(A)\times (A\times\partial D^1)$$
et
$$ \pull {Z^k(B)} \gets Z^k(B) \times (B\times\partial D^1)\to E_2^k(B)\times (B\times\partial D^1) .$$
De l'hypothèse {\bf P2} et du corollaire $\ref{PullBord}$, on déduit une équivalence d'homotopie entre $Z^{k+1}(A)$ et $Z^{k+1}(B)$ telle que le carré suivant commute à homotopie près :
$$\xymatrix{
Z^{k+1}(A)\ar[r]\ar[d]_\he	&	E_2^{k+1}(A)\ar[d]^\he	\\
Z^{k+1}(B)\ar[r]		&	E_2^{k+1}(B)	\\
}$$
Il nous reste donc à vérifier que l'espace $Z^{k+1}(A)$ satisfait les propriétés {\bf P1} et {\bf P3} de l'hypothèse de récurrence.

Par hypothèse de récurrence, $E_2^k(A)$ est la colimite homotopique de $$\FF_k (A) \gets{\FF_k (A)\times\partial D^{k-2}}\to Z^k(A)$$
donc, d'après \ref{trivcolim}, $E_2^k\times(A\times\partial D^1)$ est la colimite  homotopique de $$\FF_k (A)\times(A\times\partial D^1) \gets{\FF_k (A)\times\partial D^{k-2}}\times(A\times\partial D^1)\to Z^k(A)\times(A\times\partial D^1) .$$
Dans le diagramme suivant, le carré de gauche est donc cocartésien
$$\xymatrix{
{\FF_k (A)\times\partial D^{k-2}}\times(A\times\partial D^1)\ar[r]\ar[d]	& Z^k(A)\times(A\times\partial D^1)\ar[r]	\ar[d]&	{\pull {Z^k(A)} }\ar[d] \\
{\FF_k (A)\times(A\times\partial D^1)\ar[r]}	& E_2^k(A)\times(A\times\partial D^1)\ar[r]	&	Z^{k+1}(A)	\\
}$$
Le carré de droite étant aussi cocartésien (par définition de $Z^{k+1}(A)$), le diagramme total est cocartésien.

Notons $g$ la composée ${\FF_k (A)} \times \partial D^{k-2} \to Z^k(A)\to {\FF_k (A\times J)}$ et utilisons la pour définir l'application
$$\begin{array}{rcl}
{\FF_{k+1} (A)} \times \partial D^{k-2}\times [-k,k] &\to& {\FF_{k+1} (A\times J)} \\ (a_1,\dots,a_{k+1},z,t)&\mapsto& (g(a_1,\dots,a_k,z),(a_{k+1},t)) \\
\end{array}$$
Remarquons que cette dernière est bien définie car on a supposé que la composée 
$${\FF_k (A)\times\partial D^{k-2}}\to Z^k(A)\to {\FF_k (A\times J)}\to A^{\times k}$$
est égale à ${\FF_k (A)\times\partial D^{k-2}}\to\FF_k (A) \hookrightarrow A^{\times k}$.
La propriété universelle d'un carré cartésien nous donne une flèche hachurée dans le diagramme, commutatif à homotopie près, ci-dessous :
$$\xymatrix{
{\FF_{k+1} (A)} \times \partial D^{k-2}\times [-k,k] \ar@/^/ [rrd]\ar@{-->}[rd]\ar@/_/[d] \\
{\FF_k (A)} \times \partial D^{k-2}\ar@/_/[rd] &{\pull {Z^k(A)} }\ar[d]\ar[r]	&	{\FF_{k+1} (A\times J)}\ar[d] \\
&Z^k(A)\ar[r]				&	{\FF_k (A\times J)}	\\
}$$
Ceci montre que le carré supérieur gauche du diagramme suivant, noté $(\Xi)$, commute à homotopie près.
$$\xymatrix{
{\pull Z^k(A) }	& {\FF_k (A)\times\partial D^{k-2}}\times(A\times\partial D^1)\ar[r]\ar[l]	&	{\FF_k (A)\times(A\times\partial D^1)}	\\
{\FF_{k+1} (A) \times \partial D^{k-2}\times [-k,k]}\ar[u]\ar[d]	& {\FF_{k+1} (A) \times \partial D^{k-2}\times\partial D^1}\ar[r]\ar[l]\ar[u]\ar[d]	& {\FF_{k+1} (A) \times \partial D^{1}}\ar[u]\ar[d]	\\
{\FF_{k+1} (A) \times \partial D^{k-2}}	& {\FF_{k+1} (A) \times \partial D^{k-2}}\ar[r]\ar[l]	& {\FF_{k+1} (A)}
}$$
Les autres applications du diagramme $(\Xi)$ soit ont déjà été rencontrées auparavant, soit sont des projections sur certains facteurs d'un produit cartésien soit encore sont des inclusions naturelles. Il est facile de vérifier que le diagramme $(\Xi)$ commute à homotopie près. 
Remarquons alors que :
\begin{itemize}
	\item La colimite homotopique de la première ligne est $Z^{k+1}(A)$.
	\item La colimite homotopique de la deuxième ligne est ${\FF_{k+1} (A) \times \partial D^{k-1}}$.
	\item La colimite homotopique de la troisième ligne est $\FF_{k+1} (A)$.
	\item La colimite homotopique de la première colonne est trivialement $ \pull Z^k(A) $.
	\item La colimite homotopique de la deuxième colonne est $ \pull {\FF_k (A)\times\partial D^{k-2} } $ pour l'application $g : {\FF_k (A)\times\partial D^{k-2} }\to\FF_k(A\times J)$. (appliquer \ref{Puppe2} de manière similaire à la démonstration de \ref{ThSusp1}).
	\item La colimite homotopique de la troisième colonne est $ \pull {\FF_k (A)} $ pour l'application $\FF_k(A)\to  E^k_2(A)\to\FF_k(A\times J)$. (appliquer \ref{Puppe2} de manière similaire à la démonstration de \ref{ThSusp1}).
\end{itemize}
La colimite homotopique de notre diagramme est donc la même que celle de 
$$ \pull {Z^k(A)}  \gets \pull {\FF_k (A)\times\partial D^{k-2} } \to \pull {\FF_k (A)} $$
qui, d'après le lemme du cube \cite{Mather1976,Doeraene1998}, est égale à $\pull E_2^k(A) = E_2^{k+1}(A) $.
Cette colimite homotopique coïncide aussi avec celle de 
$$ \FF_{k+1} (A)  \gets {\FF_{k+1} (A) \times \partial D^{k-1}} \to Z^{k+1}(A) $$
ce qui établit {\bf P1}.
Notons, par ailleurs, que la composée $${\FF_{k+1} (A) \times \partial D^{k-1}}\to Z^{k+1}(A)\to {\FF_{k+1} (A\times J)}\to A^{\times k+1}$$ est bien égale à $$\FF_{k+1} (A)\times \partial D^{k-1} \to\FF_{k+1} (A)\hookrightarrow A^{\times k+1}.$$
(Il suffit de voir que la $j$-ième particule conserve sa composante en $A$ par n'importe quelle application du diagramme $(\Xi)$).
Ceci conclut notre récurrence.

D'après {\bf P1}, la cofibre de $ Z^{k+1}(A)\to E_2^{k+1}(A) $ a le type d'homotopie de la cofibre de ${\FF_{k+1} (A) \times \partial D^{k-1}}\to\FF_{k+1} (A)$ donc de $ \Sigma^{k-1}\ \FF_{k+1} (A)_+$.
L'invariance homotopique de cette suspension découle de la propriété {\bf P2} de l'hypothèse de récurrence, à savoir la commutativité à homotopie près de :
$$\xymatrix{
Z^{k+1}(A)\ar[r]\ar[d]_\he	&	E_2^{k+1}(A)\ar[d]^\he	\\
Z^{k+1}(B)\ar[r]		&	E_2^{k+1}(B)	\\
}$$
\qed

\appendix
\section{Limites Homotopiques}
\label{sec:HoLim}

Nous regroupons dans cette section quelques propriétés qui nous ont été utiles tout au long de ce travail.
Rappelons tout d'abord les définitions de produit fibré homotopique et de somme amalgamée homotopique. 
Si $H$ est une homotopie entre deux applications $f$ et $g$, nous écrivons $H : f \sim g$.

Supposons donné un diagramme commutatif à homotopie près
$$\xymatrix{
P\ar[d]_{g_1}\ar[r]^{f_1}	&	C\ar[d]^g \\
A\ar[r]_f				&B	\\
}$$
avec $H:gf_1\sim fg_1$. Un tel diagramme est un \emph{produit fibré homotopique} si pour tout espace $D$ et tout couple d'applications  $(f_2:D\to C , g_2:D\to A)$ faisant commuter à homotopie près le carré extérieur du diagramme suivant :
$$\xymatrix{
D\ar@/^/[rrd]^{f_2}\ar@{-->}[rd]|w\ar@/_/[ddr]_{g_2} \\
&P\ar[d]_{g_1}\ar[r]^{f_1}	&	C\ar[d]^g \\
&A\ar[r]_f				&B	\\
}$$
avec $G:gf_2\sim fg_2$, alors on a les propriétés suivantes :
\begin{itemize}
	\item Il existe une application $w:D\to P$ qui fait commuter la totalité du diagramme précédent à homotopie près ; c'est-à-dire que l'on a des homotopies $K:f_2\sim f_1w$ et $L:g_1w\sim g_2$ telles que $g\circ K+H\circ w + f\circ L\sim G$.
	\item L'application $w$ est unique à homotopie près ; c'est-à-dire que si on se donne une autre application $w':D\to P$ avec des homotopies  $K':f_2\sim f_1w'$ et $L':g_1w'\sim g_2$ telles que $g\circ K'+H\circ w' + f\circ L'\sim G$, alors il existe une homotopie $M:w\sim w'$ telle que $K+f_1\circ M\sim K'$ et $g_1\circ M+L'\sim L$.
\end{itemize}
On dira aussi que le carré $P-A-B-C$ est \emph{cartésien}.

En dualisant la définition précédente (i.e. en inversant le sens des flèches dans le diagramme précédent), on obtient la définition d'une \emph{somme amalgamée homotopique}.
Le diagramme sous-jacent est donné par :
$$\xymatrix{
B\ar[d]_g\ar[r]^f	&	A\ar[d]\ar@/^/[ddr] \\
C\ar[r]\ar@/_/ [rrd]		&	Z\ar@{-->}[rd]|w	\\
			&			& D	\\
}$$
On dit encore que le carré $A-B-C-Z$ est \emph{cocartésien}.

Ces notions se généralisent en \emph{limite homotopique} et \emph{colimite homotopique} \cite{Bousfield1972}. Par exemple, dans les définitions ci-dessus $P$ est la limite homotopique de $A\stackrel{f}{\rightarrow}B\stackrel{g}{\leftarrow}C$ et $Z$ est la colimite homotopique de $A\stackrel{f}{\leftarrow}B\stackrel{g}{\rightarrow}C$. 

La propriété suivante est classique et nous est souvent utile dans nos démonstrations.



\begin{PROP} [Lemme du prisme, \cite{Mather1976,Doeraene1998}] \label{lemprism} Considérons le diagramme suivant :
$$\xymatrix{
E_1\ar[rr]\ar@{-->}[dr]\ar[dd]	&	&E_3\ar[dd] 	\\
	&E_2\ar[ur]\ar[dd]	&	\\
E_1'\ar[rr]\ar[dr]	&	&E_3'	\\
	&E_2'\ar[ur]	&	\\
}$$
On suppose que ce diagramme commute à homotopie près et que $E_2-E_3-E_3'-E_2'$ est cartésien.
La définition du produit fibré homotopique nous assure l'existence d'une application $E_1\to E_2$ qui commute avec le reste du diagramme.

Alors  le carré  $E_1-E_3-E_3'-E_1'$  est cartésien si et seulement si le carré $E_1-E_2-E_2'-E_1'$ est lui aussi cartésien.
De plus, si on suppose que $E_1'\to E_2'$ est une équivalence d'homotopie alors il en est de même pour $E_1\to E_2$.
\end{PROP}


Le théorème suivant justifie la validité de notre construction. Il met en relation colimite homotopique d'applications de même fibre (\ref{Puppe}) ou de même base (\ref{Puppe2}).
\Th [V.~Puppe, \cite{Puppe1974}] \label{Puppe}
Soit $F$ un espace fixé et soit $\mathrm{Top}^F$ la catégorie des applications ayant pour fibre homotopique $F$. On considère alors un diagramme dans $\mathrm{Top}^F$ indexé par une petite catégorie $I$,
c'est-à-dire que pour tout objet $i$ de $I$, on a une application $E_i\to B_i$ ayant $F$ pour fibre homotopique et pour tout morphisme $(i,j)$ de $I$, on a une famille d'applications ${\alpha_{i,j}}$ et ${\alpha'_{i,j}}$ qui rendent le carré suivant commutatif à homotopie près :
$$\xymatrix{
E_i\ar[rr]^{\alpha_{i,j}}\ar[d]	&&	E_j\ar[d]	\\
B_i\ar[rr]^{\alpha'_{i,j}}	&&	B_j	\\
}$$
Alors l'application induite $(\textrm{hocolim}_I~ \alpha_{i,j})\to (\textrm{hocolim}_I~ \alpha'_{i,j})$ est aussi dans $\mathrm{Top}^F$.

\Cor \label{trivcolim}
Soit $F$ un espace fixé. On considère un diagramme d'espaces topologiques $(B_i)_{i\in I}$ indexés par une petite catégorie $I$. Alors
$$ \textrm{hocolim}_I (B_i\times F)=\textrm{hocolim}_I (B_i)\times F .$$


\Prop [\cite{Puppe1974}] \label{Puppe2} 
Soit $B$ un espace fixé. On considère un diagramme d'applications $\alpha_{i,j} : E_i\to E_j$ indexées par un petite catégorie $I$. On suppose que pour tout morphisme $(i,j)$ de $I$, on a le diagramme commutatif suivant :
$$
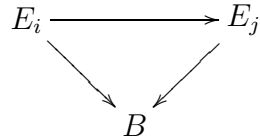
\xymatrix{ 
E_i\ar[rr]\ar[rd]		&&	E_j\ar[dl] \\
	&	B	& \\
}$$
Si l'on note $\alpha_{i,j}' : F_i\to F_j$ les applications induites au niveau des fibres homotopiques de $E_i\to B$ et $E_j\to B$, alors la fibre homotopique de $(\textrm{hocolim}_I~ \alpha_{i,j}) \to B$ a le même type d'homotopie que la colimite homotopique des $\alpha'_{i,j}$.

La propriété suivante se prête particulièrement bien à notre construction par récurrence des espaces $E^k_{\ibar}$.
Ne l'ayant pas trouvée dans la littérature, nous en détaillons la preuve.
\Prop \label{bicube} Considèrons le diagramme suivant, commutatif à homotopie près, dans lequel les applications désignées par $\he$ sont des équivalences d'homotopie.
$$
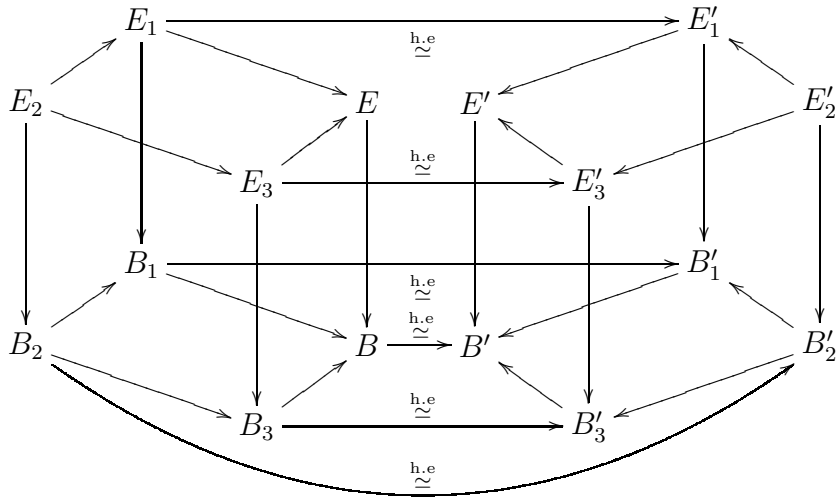
\xymatrix@R=2.5ex{
	&E_1\ar[rrrrr]_\he\ar[rrd]\ar[ddd]&&&&&E'_1\ar[lld]\ar[ddd]\\
E_2\ar[ru]\ar[rrd]\ar[ddd]	&&&E\ar[ddd]&E'\ar[ddd]&&&E'_2\ar[lu]\ar[lld]\ar[ddd]\\
	&	&E_3\ar[rrr]^\he\ar[ru]\ar[ddd]&&&E_3'\ar[lu]\ar[ddd]\\
	&B_1\ar[rrd]\ar[rrrrr]_\he&&&&&B'_1\ar[dll]	\\
B_2\ar[ru]\ar[rrd]\ar@/_11ex/[rrrrrrr]^\he	&	&	&B\ar[r]^\he&B'&&&B'_2\ar[dll]\ar[ul]\\
	&	&B_3\ar[ru]\ar[rrr]^\he&&&B_3'\ar[ul]\\}$$
Pour les cubes droit et gauche, on suppose que les faces latérales sont cartésiennes et les faces du dessus et du dessous sont cocartésiennes (ce sont des cubes du lemme du cube de M.~Mather \cite{Mather1976}). Alors, il existe des équivalences d'homotopie $E_2\to E'_2$ et $E\to E'$ qui commutent à homotopie près avec le reste du diagramme.
 
\Dem 
Appliquons le lemme du prisme (\ref{lemprism}) en choisissant pour base le triangle
$$\xymatrix{
B_2\ar[rr]\ar[dr]		&	&	B_3'	\\
		& B_2'\ar[ur]	&\\
}$$
qui commute à homotopie près.
Nous obtenons une équivalence d'homotopie $E_2\to E'_2$ qui fait apparaître $E_2$ comme le produit fibré homotopique de $E'_2\to B'_2$ le long de $B_2\to B'_2$ et fait commuter à homotopie près le carré
$$\xymatrix{
E_2\ar[r]\ar[d]	&	E'_2\ar[d]	\\
E_3\ar[r]	&	E'_3	\\
}$$

En procédant de même à partir du triangle 
$$\xymatrix{
B_2\ar[rr]\ar[dr]		&	&	B_1'	\\
		& B_2'\ar[ur]	&\\
}$$
on obtient un carré commutatif 
$$\xymatrix{
E_2\ar[r]\ar[d]	&	E'_2\ar[d]	\\
E_1\ar[r]	&	E'_1	\\
}$$
Remarquons que l'application $E_2\to E'_2$ est, à homotopie près, la même que pré\-cé\-dem\-ment car elle est obtenue par produit fibré homotopique de $E'_2\to B'_2$ le long de $B_2\to B'_2$.
En passant à la colimite, on en déduit l'équivalence d'homotopie $E\to E'$ et la commutativité du carré
$$\xymatrix{
E\ar[r]\ar[d]	&	E'\ar[d]	\\
B\ar[r]	&	B'	\\
}$$
\qed

\begin{flushleft}

\eject
\footnotesize
\bibliographystyle{abbrv}
\bibliography{FAxR}

\begin{thebibliography}{10}

\bibitem{Aouina2003}
M.~Aouina and J.~R. Klein.
\newblock On the homotopy invariance of configuration spaces.
\newblock http://arxiv.org/abs/math.AT/0310483/, 2003.

\bibitem{Benlian1967}
R.~Benlian and J.~Wagoner.
\newblock Type d'homotopie fibr\'e et r\'eduction structurale des fibr\'es
  vectoriels.
\newblock {\em C. R. Acad. Sci. Paris S\'er. A-B}, 265:A207--A209, 1967.

\bibitem{Bousfield1972}
A.~K. Bousfield and D.~M. Kan.
\newblock {\em Homotopy limits, completions and localizations}.
\newblock Springer-Verlag, Berlin, 1972.
\newblock Lecture Notes in Mathematics, Vol. 304.

\bibitem{Cohen2002}
F.~R. Cohen and S.~Gitler.
\newblock On loop spaces of configuration spaces.
\newblock {\em Trans. Amer. Math. Soc.}, 354(5):1705--1748 (electronic), 2002.

\bibitem{Doeraene1998}
J.-P. Doeraene.
\newblock Homotopy pull backs, homotopy push outs and joins.
\newblock {\em Bull. Belg. Math. Soc. Simon Stevin}, 5(1):15--37, 1998.

\bibitem{Fadell1962}
E.~Fadell and L.~Neuwirth.
\newblock Configuration spaces.
\newblock {\em Math. Scand.}, 10:111--118, 1962.

\bibitem{Greene1978}
R.~E. Greene.
\newblock Complete metrics of bounded curvature on noncompact manifolds.
\newblock {\em Arch. Math. (Basel)}, 31(1):89--95, 1978/79.

\bibitem{Klein1999}
J.~R. Klein.
\newblock Poincar\'e duality embeddings and fiberwise homotopy theory.
\newblock {\em Topology}, 38(3):597--620, 1999.

\bibitem{Klein2002}
J.~R. Klein.
\newblock Poincar\'e duality embeddings and fibrewise homotopy theory. {II}.
\newblock {\em Q. J. Math.}, 53(3):319--335, 2002.

\bibitem{Kriz1994}
I.~K{\v{r}}{\'{\i}}{\v{z}}.
\newblock On the rational homotopy type of configuration spaces.
\newblock {\em Ann. of Math. (2)}, 139(2):227--237, 1994.

\bibitem{Lambrechts2001}
P.~Lambrechts and D.~Stanley.
\newblock The rational homotopy type of the configuration space of two points.
\newblock http://gauss.math.ucl.ac.be/\~{}topalg/preprints/, 2001.

\bibitem{Levitt1995}
N.~Levitt.
\newblock Spaces of arcs and configuration spaces of manifolds.
\newblock {\em Topology}, 34(1):217--230, 1995.

\bibitem{Longoni2004}
R.~Longoni and P.~Salvatore.
\newblock Configuration spaces are not homotopy invariant.
\newblock http://arxiv.org/abs/math.AT/0401075/, 2004.

\bibitem{Mather1976}
M.~Mather.
\newblock Pull-backs in homotopy theory.
\newblock {\em Canad. J. Math.}, 28(2):225--263, 1976.

\bibitem{May1975}
J.~P. May.
\newblock Classifying spaces and fibrations.
\newblock {\em Mem. Amer. Math. Soc.}, 1(1, 155):xiii+98, 1975.

\bibitem{Puppe1974}
V.~Puppe.
\newblock A remark on ``homotopy fibrations''.
\newblock {\em Manuscripta Math.}, 12:113--120, 1974.

\bibitem{Totaro1996}
B.~Totaro.
\newblock Configuration spaces of algebraic varieties.
\newblock {\em Topology}, 35(4):1057--1067, 1996.

\end{thebibliography}

\end{flushleft}

\end{document}